\documentclass[a4paper,12pt]{article}

\usepackage{a4wide}
\usepackage{graphicx,subfigure,hhline}
\usepackage[colorlinks=true,linkcolor=blue,citecolor=red]{hyperref}
\usepackage{xcolor}
\usepackage{amssymb,amsthm,empheq,bbold}
\usepackage[inline]{enumitem}
\usepackage[ruled,vlined,linesnumbered]{algorithm2e}
\usepackage{caption}
\newtheorem{theorem}{Theorem}

\newtheorem{assumption}{Assumption}
\newtheorem{lemma}{Lemma}
\usepackage{geometry}
\usepackage[normalem]{ulem}
\usepackage{cancel} 

\usepackage{geometry} 
\allowdisplaybreaks[1]

\newcommand{\RR}{\mathbb{R}}
\definecolor{ggreen}{cmyk}{1,     0,      1,      0}
\definecolor{myred}{cmyk}{0.1, 1, 0.5, 0}
\definecolor{dblue}{rgb}{0.4, 0.29996, 0.7} 

\newcommand{\be}{\begin{equation}}
	\newcommand{\ee}{\end{equation}}

\newcommand{\fvec}{\underline{\mathbf{f}}}

\def\pa{\partial}
\def\bv{{\bf v}}

\def\bx{{\bf x}}

\allowdisplaybreaks

\begin{document}
	
	\title{
		Reaction-diffusion systems from kinetic models for 
		bacterial communities on a leaf surface}
	
	\date{}
	
	\author{\small
		Marzia Bisi{$^{1}$}, Davide Cusseddu{$^{2,3}$},  Ana Jacinta Soares{$^{3}$},  Romina Travaglini{$^{4,1,^*}$}\\[1em]
		{\footnotesize $^1$}{\small\it Dept. of Mathematical, Physical and Computer Sciences,} \\ 
		{\small \it University of Parma, Parco Area delle Scienze 53/A, 43124, Parma, Italy}\\{\footnotesize marzia.bisi@unipr.it}
		\\[0.4em]
		{\footnotesize $^2$}{\small\it Department of Mathematical Sciences ‘‘G. L. Lagrange’’, Politecnico di Torino,}\\ 
		{\small \it Corso Duca degli Abruzzi, 24, 10129 Torino,  Italy}
		\\{\footnotesize davide.cusseddu@polito.it}
		\\[0.4em]
		{\footnotesize $^3$}{\small\it Centre of Mathematics of the University of Minho,}\\ 
		{\small \it  Campus de Gualtar, 4710-057 Braga, Portugal}
		\\{\footnotesize ajsoares@math.uminho.pt}
		\\[0.4em]
		{\footnotesize $^4$}{\small \it INDAM – Istituto Nazionale di Alta Matematica ‘‘F. Severi’’,}\\ 
		{\small \it Piazzale Aldo Moro 5, 00185, Roma, Italy}
		\\[0.4em]
		$^*${\footnotesize corresponding author }\\{\footnotesize romina.travaglini@unipr.it}
	}		
	
	\maketitle
	
	
	\begin{abstract}
		Many mathematical models for biological phenomena, such as the spread of diseases, are based on reaction-diffusion equations for densities of interacting cell populations. 
		We present a consistent derivation of reaction-diffusion equations from systems of suitably rescaled 
		{kinetic equations} for distribution functions of cell populations interacting in a host medium.
		We show at first that the classical diffusive limit of kinetic equations leads to linear diffusion terms only. 
		Then, we show possible strategies in order to obtain, from the kinetic level, macroscopic systems 
		with nonlinear diffusion and also with cross-diffusion effects. 
		The derivation from a kinetic description has the advantage of relating reaction and diffusion coefficients 
		to the microscopic parameters of the interactions. 
		We present an application of our approach to the study of the evolution of different bacterial populations on a leaf surface. 
		Turing instability properties of the relevant macroscopic systems are investigated by analytical methods and numerical tools, 	with particular emphasis on pattern formation for varying parameters in two-dimensional space domains.
	\end{abstract}
	
	
	\smallskip
	
	
	\smallskip
	
	
	
	\section{Introduction}
	\label{sec:int}
	
	The description of biological phenomena by means of mathematical models can be performed through different approaches. 
	We focus on the need to describe complex systems, composed of many heterogeneous living individuals, interacting stochastically within themselves and with the external environment, at spatial scales considerably smaller than at the observable level. One of the most suitable tools to perform such a description is the kinetic theory of active particles \cite{survey}. This approach derives from the classical kinetic theory of inert matter, whose key element is the Boltzmann equation. It extends the concept of interacting entities from binary, short-range collisions between molecules to non-local, multiple interplays of living individuals. Such interactions are described by systems of integro-differential equations, as proposed since early works as \cite{bellomo1994dynamics}. The kinetic approach resulted in being useful in describing a wide number of physical problems, ranging from medical studies \cite{della2022mathematical,ramos2019kinetic} to socio-economics \cite{bertotti2023modelling, cordier2005kinetic}; for further references, we address the reader to \cite{survey}.
	The kinetic theory of active particles is based on the fact that each entity/population involved is described by a distribution function, usually depending on time, space, velocity, and a further variable (activity) representing the particular state of microscopic interacting agents (typically cells or individuals) with respect to a specific characteristic.
	
	A further powerful feature of kinetic theory is the possibility of describing different types of interactions at multiple spatial or temporal scales. This allows, in particular, to obtain, through proper diffusive limits, partial-differential equations of reaction-diffusion type for observable quantities, such as macroscopic densities of constituents.
	Additionally, such equations allow us to investigate how the microscopic dynamics affect the global behavior at the macroscopic level.
	Some examples of this procedure may be found in the frame of classical Boltzmann theory of gas dynamics \cite{bisi2006reactive,bisi2022reaction,lachowicz2002microscopic}, but also in the kinetic theory of active particles describing cells and tissues (see \cite{burini2019multiscale} and references therein).
	
	Models cited above for cellular dynamics have been refined in order to obtain more complex diffusive terms, like the ones accounting for chemotaxis \cite{alt1980biased,othmer2000diffusion,othmer2002diffusion}, and applied to medical issues like the study of cancer \cite{conte2023mathematical} or multiple sclerosis \cite{oliveira2024reaction}.	
	In all of these models, though, diffusive terms are derived from the assumption that the dominant processes are the interactions whose result is a change in the velocity of the cell. We propose a new approach in which we consider a certain number of cell populations interacting with a host medium (inspired by kinetic models for gases diffusing in the atmosphere \cite{bisi2022reaction}). In this context, we suppose that the interactions of cells with the host are the dominant process and that they are conservative, in the sense that the outcome is a change in the cellular activity or the cell direction, but not in the number of cells. Then we take into account other types of phenomena that induce growth or decay of the cell populations, or influence the movement of cells, that occur at slower time scales.
	These assumptions, along with specific hypotheses on the velocity of cells, allow us to derive a reaction-diffusion system potentially including cross-diffusion terms for macroscopic densities of populations involved.
	
	{ Our approach builds upon and complements several recent multiscale modeling efforts that derive cross-diffusion systems from kinetic or microscopic descriptions. For instance, authors of \cite{bendahmane2024mathematical} derived a nonlinear predator–prey system with cross-diffusion and fluid interaction using multiscale techniques, while in \cite{zagour2019modeling}, kinetic and macroscopic approaches are applied to complex systems in biology and traffic flow. Moreover, in \cite{zhigun2022flux,dietrich2022multiscale} rigorous upscaling methods are derived, leading to flux-limited macroscopic PDEs for cancer invasion, incorporating receptor-based dynamics. Compared to these works, our kinetic-to-macroscopic derivation focuses on a specific ecological context and emphasizes the derivation of reaction-cross-diffusion systems from a mesoscopic description with a clear biological interpretation of the coefficients. While the resulting macroscopic equations are simpler and our upscaling approach is less general than the frameworks presented in the cited works, our aim is to provide a tractable kinetic framework that directly links microscopic interaction rules to pattern formation in a concrete biological setting.}
	
	At the macroscopic level, biological phenomena of main interest are those involving the formation of patterns. These can be found, indeed, in morphogenesis, chemistry, or landscape. The mathematical description of such dynamics may be obtained by means of Turing instability analysis of a reaction-diffusion system \cite{turing1990chemical}, occurring when a spatially homogeneous steady state turns into symmetry-breaking structures due to the presence of diffusive terms.
	
	{It is important to emphasize that the model developed here is primarily conceptual. Its main purpose is to show how macroscopic reaction-cross-diffusion systems can be systematically derived from a kinetic framework. At the same time, we propose a possible biological application, intended to illustrate how the methodology may be adapted to a concrete setting.} A particular biological process that can be described using the procedure outlined in this work is the aggregation of bacterial strains on a leaf surface. Several studies show the tendency of bacteria to aggregate in biofilms
	\cite{monier2004frequency,monier2003differential}, predominantly in the areas of the leaf where water and nutrients are more prevalent, like trichomes (secretion organs of the leaf), veins, and epidermal cell grooves \cite{brewer1991functional,morris1997methods}. Moreover, there is a wide biological literature concerning the interaction of two different bacterial populations on a leaf surface that may influence the aggregation (see \cite{schlechter2019driving} and references therein). Microbial interactions may be classified into cooperation and competition. Cooperation denotes interactions where at least one strain benefits without causing harm to others. Conversely, competitive relationships involve detrimental effects on at least one population, stemming from interference or exploitation competition. These different interplays may lead to various spatial organizations of bacterial populations, like co-aggregation,
	segregation, or random distributions \cite{schlechter2019driving}. 
	Further findings suggest that bacterial colonizers on leaves interact with their environment across various spatial scales. Interactions among bacteria tend to occur predominantly at small spatial scales, contrasting with those between bacteria and leaf surface structures, which extend noticeably beyond typical microscopic dimensions \cite{esser2015spatial}.
	
	From a mathematical modeling point of view, the classical macroscopic models describing the dynamics between two species, such as the Lotka-Volterra model \cite{wangersky1978lotka}, have been variously extended. In \cite{ma1996mathematical}, competitive populations’ behaviors have been included, while works like \cite{abbas2010existence} have also considered the effect of substances produced by bacterial interactions, which may positively affect the growth of the species. This chemical substance has been modeled in \cite{mu2023hopf}, along with the inclusion of diffusive terms, obtaining a reaction-diffusion system whose pattern formation has been analyzed. 
	
	To our knowledge, only a few studies have been devoted to modeling this particular biological problem.
	Some of these employ individual-based models \cite{van2013explaining}, {which provides a unified upscaling framework for bacterial dispersion and can handle both chemotaxis and haptotaxis,} 
	while others rely on stochastic models \cite{perez2012stochastic}.
	The novelty of the present work lies in the development of a proper kinetic framework specifically adapted to describe two interacting bacterial populations on a leaf surface and from which a macroscopic system can be derived.
	From the mathematical point of view, even if the asymptotic procedure leading from kinetic equations to reaction--diffusion systems with linear self--diffusion is well established in the literature \cite{bisi2006reactive}, the method allowing to rigorously obtain nonlinear and cross--diffusion terms is still under investigation, and this paper represents a step in this direction.
	
	The paper is organized as follows. 
	In Section \ref{SecKin}, a general kinetic setting for a certain number of cell populations interacting in a host medium (host tissue) is outlined. In Section \ref{SecDiffLim}, a reaction-diffusion system is derived from the proposed kinetic model under suitable scaling assumptions. Then, in Section \ref{SecDiffLim2}, an analogous procedure is outlined, including other operators in the kinetic equations leading to cross-diffusion at the macroscopic level.
	In Section \ref{SecApp}, the general strategy leading from kinetic to reaction-diffusion systems is applied to the case of microbial populations on a leaf surface, and Turing instability analysis of the obtained macroscopic equations is performed. 
	In Section \ref{SecSim} some numerical simulations, 
	performed by using the method sketched in Appendix \ref{Appendix:numerical_method}, are provided to validate
	theoretical results.
	Finally, Section \ref{SecConc} contains some final observations and future perspectives.
	
	
	\section{Kinetic description}
	\label{SecKin}	
	
	
	We propose a kinetic model for $N$ cellular populations $C_1$, $C_2\,\ldots\,C_N$ interacting among themselves and diffusing in a much denser cellular tissue, the host medium $H$, and spreading on a spatial domain 
	${\Gamma_{\bx}}{\subset\RR^n}$, where $n$ may take the values $1, 2, 3$. 
	A possible biological application of this frame will be shown in Section \ref{SecApp}. 
	We describe each population $C_i$ by means of a distribution function $f_i(t,\bx,\bv,u)$ depending on time $t\in[0,+\infty)$, position $\bx \in {\Gamma_{\bx}}$, cellular velocity $\bv \in \Gamma_{\bv}\subset \mathbb{R}^n$, 
	and on an activity variable $u$ belonging to a set $\Sigma$ which is assumed symmetric with respect to $u=0$.
	
	As a novel feature of this work, we assume that the speed of each cell depends on its position, activity, and time, 
	being the dependence law specific for each population. There is no constraint, instead, on the cell direction. This modeling choice is motivated by our objectives, as further explained in Section \ref{SecApp}.
	Thus, by decomposing the velocity variable as $\bv=v\,\hat\bv$, with $v$ being the speed and $\hat\bv$ the direction, we may write {it} as $\bv=\bv(t, \bx,\,u)=\hat\bv\,c_i(t,\,\bx,\,u)$, 
	with $\hat\bv\in\mathbb{S}^{n-1}$, and  $v=c_i(t,\,\bx,\,u)$ representing the cellular speed of population $C_i$,  {with $c_i:[0,+\infty)\times{\Gamma_{\bx}}\times\Sigma\rightarrow[0,\bar c]$, being $\bar c$ the maximal speed of cells, thus $ \Gamma_{\bv}=\bar c\,\mathbb S^{n-1}$}. 
	Distribution functions can be thus expressed as $f_i(t,\bx,\hat\bv,u),\, i=1,\ldots,N$.
	
	For the host medium, instead, we suppose that the cells of this population exist in a huge quantity, so that its distribution $f_H$ is uniform in time, space, and velocity,
	and it just depends on cellular activity $u\in\Sigma^H$, with $\Sigma^H$ being symmetric with respect to $u=0$. Indeed, since the medium $H$ is much denser, we assume \( f_H \) to be a given function, determined at the macroscopic scale. In particular, it is considered isotropic with respect to \( u \), and not influenced by interactions with the cell populations \( C_1, \dots, C_N \). This reflects the idea that the background medium remains in a stationary, homogeneous state that can be externally measured or estimated. This is a usual assumption in the kinetic description of molecules diffusing in a background medium, see for example \cite{bardos2016simultaneous, bisi2006reactive, canizo2018rate}.
	
	Densities of cellular populations may be recovered as appropriate moments of the distribution functions. Specifically,
	\begin{equation} \label{ni}
		n_i(t,\bx)=\int_{\mathbb{S}^{n-1}}\int_{\Sigma}f_i(t,\bx,\hat\bv,u)\,du\,d\hat\bv, \quad
		i=1,\ldots,N,
	\end{equation}
	\noindent
	provides the total density of population $C_i$ at time $t$ and position $\bx$.
	Analogously,
	\begin{equation}  n_H=\int_{\Sigma^H} f_H(u)\,du 
	\end{equation}
	yields the total density of the host tissue, for which we additionally suppose that the mean activity is zero  (isotropy property), i.e.
	\begin{equation}
		\int_{\Sigma^H} u\,f_H(u)\,du=0,
	\end{equation}
	since the host medium has a huge quantity of cells, but it can be considered in an equilibrium (steady) state in the absence of external populations acting on it; therefore, the mean global effect in terms of cellular activity is almost imperceptible.

	The evolution of each distribution function is described by an integro-differential kinetic equation given by
	{
		\begin{equation}
			\label{EqKin}
			\frac{\pa f_i}{\pa t}+\nabla_{\bx}\cdot \left(c_i\,\hat\bv f_i\right){+\,\nabla_{\bv}\cdot\left(f_i\,\mathbf S_i \right)}=\mathcal G^H_i[f_i,f_H]+\mathcal H_i[\fvec] ,{\quad i=1,\ldots,N,}
		\end{equation}
		being $\fvec$ the function vector $(f_1,\ldots,f_N)$ for populations $C_1, \dots, C_N$, respectively and $f_H$ the distribution function of the host medium introduced above. Moreover, due to the dependence of the cell velocity $\bv$ on time $t$, we include in the kinetic equation a transport term with respect to $\bv$, accounting for the acceleration of cells, given by
		\be\label{Esse}
		\mathbf S_i =\frac{\partial c_i}{\partial t}\,\hat\bv.
		\ee}
	The terms on the left-hand side of equation \eqref{EqKin} describe the free motion of cells in the absence of interactions, 
	while those on the right-hand side describe the interactive processes among cells.
	Specifically,
	$\mathcal G_i^H$ is an integral operator that takes into account the fact that the activity and the direction of each cell may change 
	through interactions with the host medium.
	Operator $\mathcal H_i[\fvec]$ accounts for the effects on population $C_i$ 
	due to natural birth and death processes, and to interactions among populations $C_1, \dots, C_N$.
	
	{
		Before introducing the specific form of the interaction terms, we recall here some fundamental differences between physical systems and living organisms. 
		In fact, unlike the Boltzmann equation, which assumes binary, short-range interactions in a rarefied medium, interactions in biological systems are often non-local, perception-driven, and involve multiple individuals. Moreover, conservation laws typical of mechanical systems generally do not hold in living systems. For these reasons, we define in the next subsections the detailed expressions of the interaction operators, along with their properties, taking into account the specific nature of the biological interactions under consideration.}
	
	\subsection{Conservative terms}\label{ConsDyn}
	From this point onward, in the definition of operators, we omit the dependence on $(t, \bx)$ for a lighter notation. The operator $\mathcal G_i^H$ accounting for the conservative interactions is given by
	\begin{equation}
		\label{ConsOp}
		\begin{aligned}
			\mathcal G_i^H[f_i,f_H]&(\hat\bv,u)
			= \iiint\limits_{\mathbb S^{n-1}\Sigma^H\times\Sigma}\  \Big[\eta_{i}^H(\hat\bv',u',u_*)\beta_i^H(\hat\bv,u;\hat\bv', u',u_*)\, f_i(\hat\bv',u') \\[0mm]
			&    -\eta_{i}^H(\hat\bv,u,u_*)\beta_i^H(\hat\bv',u';\hat\bv,\, u,u_*)\,f_i(\hat\bv,u) \Big]\,f_H(u_*)\,du'\,du_*\,d\hat\bv',  {\mbox{ for all }   (t,\bx)  \in [0,+\infty]  \times {\Gamma_{\bx}},}
		\end{aligned}
	\end{equation}
	{where  {in the loss term (second line)
			\begin{itemize}
				\item[--] $\hat\bv\in\mathbb{S}^{n-1}$ and $u\in \Sigma$ are the velocity direction and the activity, respectively, of cellular population $C_i$ before the interaction with the host;
				\item[--] $u_*\in\Sigma^H$ is the activity of the host environment cell participating in the interaction; 
				\item[--] $\hat\bv'\in\mathbb{S}^{n-1}$ and $u'\in \Sigma$ are the velocity direction and the activity, respectively, of cellular population $C_i$ after the interaction with the host.
			\end{itemize}
			In the gain term, pre- and post-interaction variables are exchanged, as usual in kinetic operators.}} 
	Terms $\eta_i^H(\hat\bv,u,u_*)\geq 0$ are the interaction frequencies between 
	a cell of the population $C_i$ having activity $u$ and velocity directed along $\hat\bv$ and a host cell having activity $u_*$, 
	whereas terms $\beta_i^H(\hat\bv,u; \hat\bv',u',u_*) $ and $\beta_i^H(\hat\bv',u';\hat\bv,\, u,u_*)$ are the transition probabilities for a cell $C_i$ 
	to pass from activity $u'$ and velocity $\hat\bv'$ to activity $u$ and velocity $\hat\bv$ or vice-versa after interaction with a host cell having activity $u_*$. For simplicity, we suppose that these functions only involve the activity (i.e.\ the cell changes its velocity with uniform probability as a result of the interaction with the host cell),
	and that $\beta_i^H$ fulfills
	\begin{equation}
		\label{prob-beta}
		\iint\limits_{\mathbb S^{n-1}\Sigma}\beta_i^H(u; u',u_*)\,du\,d\hat\bv = 1, {\quad \mbox { for all }  (u',u_*) \in \Sigma \times\Sigma^H}.
	\end{equation}
	We also observe that the interaction process of cells with the host medium is conservative, i.e.
	\begin{equation}
		\iint\limits_{\mathbb{S}^{n-1} \Sigma} \mathcal G_i^H[f_i,f_H](t,\bx,\hat\bv,u)\, du\,d\hat\bv =0, {\quad \mbox { for all }  {(t,\bx)}  \in [0,+\infty]  \times {\Gamma_{\bx}},}
	\end{equation}
	meaning that there is no direct overall proliferation or destruction of cells for each population $C_i$ resulting from their interactions with the host medium.
	\subsubsection*{{Properties of the conservative operators}}
	{In view of the asymptotic analysis developed in Sections \ref{SecDiffLim} and \ref{SecDiffLim2}, some properties and rigorous results regarding the conservative operators are needed. To this aim, we start by stating the following assumption regarding the existence of a proper detailed balance for the operators $\mathcal G^H_i$.}
	\begin{assumption}
		\label{Assu}
		{For each $i=1,\ldots,N$, let $\mathcal G_i^H[f_i,f_H]$} be the conservative operator defined in equation \eqref{ConsOp}.
		Then, there exists a {distribution} function $M_i>0$ defined on {$\mathbb{S}^{n-1}\times\Sigma$, uniform in $\hat\bv$,}
		and independent of $\bx$ and $t$, such that
		\begin{eqnarray}
			\lefteqn{\int\limits_{\Sigma^H} 
				\Big[ \eta_{i}^H(w,u_*)\beta_i^H(u; w,u_*) M_i(w)} \\
			& &  \hspace*{1cm}
			- \eta_{i}^H(u,u_*)\beta_i^H(w; u,u_*) M_i(u) \big] f_H(u_*) du_* 
			= 0 {\quad \mbox { for all }  (w,u) \in \Sigma \times\Sigma}.
			\nonumber
		\end{eqnarray}
		{Such distribution is} normalized {and its first moment in the activity variable vanishes},
		that is
		\begin{equation}\label{PropM}
			\iint\limits_{\mathbb{S}^{n-1} \Sigma} M_i(u)\,d\hat\bv\,du=1,\qquad \int_{\Sigma} u\,M_i(u)\,du=0.\end{equation}
		Moreover, there exists a constant $\gamma >0$ such that the following bound {condition} holds,
		\begin{equation}
			\label{ineq-Assu}
			\int_{\Sigma^H}\eta_{i}^H(w,u_*)\beta_i^H(u; w,u_*)\,f_H(u_*)\,du_* \geq \gamma\, M_i(u) ,
			\quad  {\mbox { for all }  (w,u) \in \Sigma \times\Sigma}.
		\end{equation}
	\end{assumption}
	
	From an applied point of view, with Assumption \ref{Assu}, we consider, for each of the $N$ populations, 
	the existence of a configuration $M_i$ in which the population remains in a state of activity equilibrium 
	with respect to the interactions with the host medium. In particular, the mean activity of each population
	in this configuration is assumed to be zero, as is the one for the host medium.
	
	\medskip
	The previous assumptions allow us to state a key result in the following lemma, where a solvability condition for 
	{integral equations involving} the operators $\mathcal G^H_i$ is provided to explicitly determine the terms of the asymptotic expansion of distribution functions performed in Sections \ref{SecDiffLim} and \ref{SecDiffLim2}.
	In the proof of Lemma \ref{Lem}, we use the following 
	Lax-Milgram theorem { (see, e.g., \cite{brezis2010functional})}.
	
	\begin{theorem} 
		Consider a bilinear form $B: \mathbb H \times \mathbb H \rightarrow \mathbb{R}$ on the Hilbert space $\mathbb H$. Assume $B$ continuous,
		i.e. there exists $C \geq 0$ such that $|B(x,y)| \leq C \|x \| \| y\|$  for all  $x, y \in \mathbb H$ , and coercive,
		i.e. there exists  $\gamma>0$ such that $B(x,x) \geq \gamma \| x \|^2$  for all  $x \in \mathbb H$. 
		Then, given any $w \in \mathbb H $, 
		there exists a unique element $x \in \mathbb H$ such that $B(u, x) = \langle u, w \rangle$  for all  $u \in \mathbb H$.
		\label{th:LM}
	\end{theorem}
	
	\noindent
	Now, we can state the following result.
	
	\begin{lemma}
		\label{Lem}
		{For each  $i=1,\ldots,N$,  let $\mathcal G_i^H[\cdot,f_H]$ be the operators  defined in \eqref{ConsOp} and let  Assumption \ref{Assu} hold. Let $g_i\in L^2\left(\mathbb S^{n-1}\times\Sigma,\dfrac{du\,d\hat\bv}{M_i}\right)$ satisfy the solvability condition: 		
			$$ \iint\limits_{\mathbb S^{n-1}\Sigma} g_i(\hat\bv,\,u)\,du\,d\hat\bv=0. $$ Then{, for any fixed $f_H\in L^2\left(\Sigma^H\right)$,} there exists a unique $h_i\in L^2\left(\mathbb S^{n-1}\times\Sigma,\dfrac{du\,d\hat\bv}{M_i}\right)$ such that 		
			\begin{equation}
				\mathcal G_i^H[h_i,f_H]=g_i, \quad {\mbox { for all } \hat\bv\in \mathbb S^{n-1}}, u\in \Sigma,
			\end{equation} satisfying $\displaystyle\iint\limits_{\mathbb S^{n-1}\Sigma} h_i(\hat\bv,u)\,du\,d\hat\bv=0$.}
	\end{lemma}
	\proof
	The assumption $\iint\limits_{\mathbb S^{n-1}\Sigma} g_i(\hat\bv,\,u)\,du\,d\hat\bv=0$ is necessary for the solvability of equation $\mathcal G_i^H[h_i,f_H]=g_i$, 
	since the linear operator $\mathcal G_i^H[h_i,f_H]$ guarantees conservation of the number of cells of the population $C_i$.
	
	To prove that this assumption is also sufficient, for any $i=1, \dots, N$, let $M_i(u)$ be a function satisfying Assumption \ref{Assu}. {Let} us consider the {following} term,
	{where}, for brevity, we omit the dependence on $u_*$ 
	of quantities involved and the dependence on pre-interaction activities of $\beta_i^H$.
	{We have{, for any fixed $f_H\in L^2\left(\Sigma^H\right)$,}}
	\begin{equation}\label{16}
		\begin{array}{cl}
			& \displaystyle \iint\limits_{\mathbb S^{n-1}\Sigma} \mathcal G_i^H[h_i,f_H](\hat\bv,\,u)\, \frac{h_i(\hat\bv,\,u)}{M_i(u)}\, du\,d\hat\bv  \\[4mm]
			& \quad =\displaystyle \iint\limits_{\mathbb S^{n-1}\times\mathbb S^{n-1}}\,\iiint\limits_ {\Sigma \times \Sigma\times\Sigma^H} \left[\eta_{i}^H(w)\beta_i^H(u)\, h_i(\hat\bv',\,w)-\eta_{i}^H(u)\beta_i^H(w)\,h_i(\hat\bv,\,u)\right]
			\\[4mm] &\qquad\qquad\qquad \times
			f_H(u_*)\, \dfrac{h_i(\hat\bv,\,u)}{M_i(u)}\, \,du_*\,\,dw \,du\,d\hat\bv'\,d\hat\bv \\[4mm]
			&\quad = \displaystyle \iint\limits_{\mathbb S^{n-1}\times\mathbb S^{n-1}} \iiint\limits_{\Sigma \times \Sigma\times\Sigma^H}  \!
			\left[\eta_{i}^H(w)\beta_i^H(u) M_i(w)\,\frac{h_i(\hat\bv',\,w)}{M_i(w)}-\eta_{i}^H(u)\beta_i^H(w)\,M_i(u)\,\frac{h_i(\hat\bv,\,u)}{M_i(u)}\right] \\[4mm]
			&\qquad\qquad\qquad \times  f_H(u_*)\, \dfrac{h_i(\hat\bv,\,u)}{M_i(u)}\,\,du_*\,\,dw \,du\,d\hat\bv'\,d\hat\bv\\[4mm]
			& \quad= \displaystyle \iint\limits_{\mathbb S^{n-1}\times\mathbb S^{n-1}} \iiint\limits_{\Sigma \times \Sigma\times\Sigma^H} \!
			\left[\eta_{i}^H(u)\beta_i^H(w) M_i(u)\,\dfrac{h_i(\hat\bv,\,u)}{M_i(u)}-\eta_{i}^H(w)\beta_i^H(u)\,M_i(w)\,\frac{h_i(\hat\bv',\,w)}{M_i(w)}\right] \\[4mm]
			&\qquad\qquad\qquad \times f_H(u_*)\, \dfrac{h_i(\hat\bv',\,w)}{M_i(w)}\, \,du_*\,\,dw \,du\,d\hat\bv'\,d\hat\bv,
		\end{array}
	\end{equation}
	where last line has been obtained by exchanging $\left(\hat\bv,u\right) \leftrightarrow \left(\hat\bv',w\right)$. 
	By summing last two lines of formula (\ref{16}) and recalling Assumption \ref{Assu}, we get
	\begin{eqnarray}
		\lefteqn{\iint\limits_{\mathbb S^{n-1}\Sigma} \mathcal G_i^H[h_i,f_H](\hat\bv,\,u) \frac{h_i(\hat\bv,\,u)}{M_i(u)} du\,d\hat\bv
			= \frac12 \! \iint\limits_{\mathbb S^{n-1}\times\mathbb S^{n-1}}\,\iiint\limits_ {\Sigma \times \Sigma\times\Sigma^H} 
			\eta_{i}^H(u)\beta_i^H(w) M_i(u)} 
		\label{17}  \\
		& & \times \left[ 2\, \frac{h_i(\hat\bv,\,u) h_i(\hat\bv',\,w)}{M_i(u) M_i(w)} \!-\! 
		\left( \frac{h_i(\hat\bv,\,u)}{M_i(u)} \right)^{\!\!2} 
		\!-\! \left( \frac{h_i(\hat\bv',\,w)}{M_i(w)} \right)^{\!\!2} \right] 
		f_H(u_*)\,du_*\,\,dw \,du\,d\hat\bv'\,d\hat\bv.
		\nonumber
	\end{eqnarray}
	Therefore,
	\begin{align} \label{ineq-sign}
		\iint\limits_{\mathbb S^{n-1}\Sigma} \mathcal G_i^H[h_i,f_H](\hat\bv,\,u)\, &\frac{h_i(\hat\bv,\,u)}{M_i(u)}\, du\,d\hat\bv = -\, \frac12 \iint\limits_{\mathbb S^{n-1}\times\mathbb S^{n-1}}\,\iiint\limits_ {\Sigma \times \Sigma\times\Sigma^H} \eta_{i}^H(u)\beta_i^H(w) M_i(u)
		\\[0.1cm]
		&  \times  
		\nonumber
		\left( \frac{h_i(\hat\bv,\,u)}{M_i(u)} - \frac{h_i(\hat\bv',\,w)}{M_i(w)} \right)^2 f_H(u_*)\, du_*\,\,dw \,du\,d\hat\bv'\,d\hat\bv.
	\end{align}
	Owing to the inequality (\ref{ineq-Assu}) of Assumption \ref{Assu}, we note that
	\begin{equation}
		\label{boundedness-inequality}
		\begin{aligned}
			\displaystyle - \! \iint\limits_{\mathbb S^{n-1}\Sigma}& \mathcal G_i^H[h_i,f_H](\hat\bv,\,u)\, \frac{h_i(\hat\bv,\,u)}{M_i(u)}\, du\,d\hat\bv\\
			& \geq \frac{\gamma}{2} \iiiint\limits_{\mathbb S^{n-1}\times\mathbb S^{n-1} \Sigma \times \Sigma} \!
			M_i(u) M_i(w) \! \left( \frac{h_i(\hat\bv,\,u)}{M_i(u)} \!-\! \frac{h_i(\hat\bv',\,w)}{M_i(w)} \right)^{\!2} \,dw \,du\,d\hat\bv'\,d\hat\bv
			\\
			&\geq \gamma\, \iint\limits_{\mathbb S^{n-1}\Sigma} \frac{h_i^2(\hat\bv,\,u)}{M_i(u)}\, du\,d\hat\bv,
		\end{aligned}
	\end{equation}
	where use has been made of the {normalization} $\displaystyle \iint\limits_{\mathbb{S}^{n-1} \Sigma} M_i(u)\,du\,d\hat\bv= 1$ and of the
	constraint \linebreak $\displaystyle \iint\limits_{\mathbb S^{n-1}\Sigma} h_i(\hat\bv,\,u)\,du\,d\hat\bv=0$. 
	
	Now, the existence and uniqueness of a weak solution to the equation $\mathcal G_i^H[h_i,f_H]=g_i$ 
	is provided by {the} Lax-Milgram theorem, see Theorem \ref{th:LM} above.
	In our case, for any fixed $i =1, \dots, N$, we set $\mathbb H = L^2\left(\mathbb S^{n-1}\times\Sigma,\frac{du\,d\hat\bv}{M_i}\right)$ and we consider the bilinear form
	$$
	B(h,k) = - \iint\limits_{\mathbb S^{n-1}\Sigma} \mathcal G_i^H[h,f_H](\hat\bv,\,u)\, \frac{k(\hat\bv,\,u)}{M_i(u)}\, du\,d\hat\bv\,,{ \quad h,k\in\mathbb H.}
	$$ 
	{The c}ontinuity of the operator $B(h,k)$ is straightforward, 
	coerciveness in the weighted $L^2$ space 
	follows directly from condition \eqref{boundedness-inequality}, since{,  for all  $h\in\mathbb H$,}
	$$
	B(h,h) = - \iint\limits_{\mathbb S^{n-1}\Sigma} \mathcal G_i^H[h,f_H](\hat\bv,\,u)\, \frac{h(\hat\bv,\,u)}{M_i(u)}\, du\,d\hat\bv \geq \gamma \iint\limits_{\mathbb S^{n-1}\Sigma} \frac{h^2(\hat\bv,\,u)}{M_i(u)}\, du\,d\hat\bv= \gamma \|h \|^2\,.
	$$
	Setting $w_i=-g_i$, we note that, by the Lax-Milgram theorem, there exists a unique solution 
	$h_i \in L^2\left(\mathbb S^{n-1}\times\Sigma,\frac{du\,d\hat\bv}{M_i}\right)$ of the equation $B(h_i, k) = \langle w_i, k \rangle$, 
	for any $k \in L^2\left(\mathbb S^{n-1}\times\Sigma,\frac{du\,d\hat\bv}{M_i}\right)$, i.e. a unique solution of
	$$
	\iint\limits_{\mathbb S^{n-1}\Sigma} \mathcal G_i^H[h_i,f_H](\hat\bv,\,u)\, \frac{k(\hat\bv,\,u)}{M_i(u)}\, du\,d\hat\bv
	= - \int_\Sigma \frac{w_i(\bv,\,u) k(\hat\bv,\,u)}{M_i(u)}\, du\,d\hat\bv,
	$$ for any $ \displaystyle k \in L^2\left(\mathbb S^{n-1}\times\Sigma,\frac{du\,d\hat\bv}{M_i} \right) .
	$
	We conclude that such a unique solution $h_i$ is a weak solution to the equation $\mathcal G_i^H[h_i,f_H]=g_i$,
	and the proof is then complete.
	\endproof

	\subsection{Non-conservative {terms}}
	
	The interaction term $\mathcal H_i$ describing the non-conservative processes in the kinetic equation (\ref{EqKin}) may be cast as
	\begin{equation}
		\mathcal H_i[\fvec] 
		= \mathcal J_i[f_i] + \sum_{j=1}^N\mathcal N_{ij}[f_i,f_j]+\sum_{\substack{j,k=1 \\ j,k\ne i}}^N\mathcal Q_{jk}^i[f_j,f_k],
		\label{eq:ops}
	\end{equation}
	where the operators describe different processes.
	Operator $\mathcal J_i$ accounts for the natural reproduction or decay of population $C_i$ and is expressed by
	\begin{equation}
		\mathcal J_i[f_i](\hat\bv,u) =\left[\vartheta_i(\hat\bv,u)-\tau_i(\hat\bv,u)\right]f_i(\hat\bv,u),
		{\quad \mbox { for all }  {(t,\bx)}  \in [0,+\infty]  \times {\Gamma_{\bx}},}
		\label{eq:Ji}
	\end{equation}
	with $\vartheta_i$ and $\tau_i$ being the {growth} and death rates, respectively. {We remark that we do not consider a direct effect of the host medium on the growth or decay of the population size. However, these phenomena are not uncorrelated. 
		Indeed,  while birth and death rates only depend on the activity $u$ and velocity direction $\hat\bv$, interactions with the host medium result in a change in the population activity and velocity direction. 
		
		Operators $\mathcal N_{ij}$ in \eqref {eq:ops}, instead, are integral operators related to interactions among cells 
		of the reference population $C_i$ 
		and cells of only another population $C_j$, including the case $j=i$. 
		We allow these interactions to be non-conservative, namely, proliferative or destructive 
		for population $C_i$.
		Thus, we take operators in a more general form with respect to the conservative ones given in (\ref{ConsOp}), 
		that is
		\begin{equation}
			\label{coll_ope}
			\begin{aligned}
				\mathcal N_{ij}[f_i,f_j]  (\hat\bv,u) 
				= \iint\limits_{\mathbb S^{n-1}\times\mathbb S^{n-1}} & \iint\limits_{\Sigma\times\Sigma}
				\mu_{ij}(\hat\bv',\hat\bv_*, u',u_*)\varphi_{ij}(\hat\bv, u; \hat\bv',\hat\bv_*, u',u_*) \\
				& \qquad \qquad \times f_i(\hat\bv',u')f_j(\hat\bv_*,u_*)du_*\,du'\,d\hat\bv_*\,d\hat\bv' \\[2mm]	      
				& - f_i(\hat\bv,u)\iint\limits_{\mathbb{S}^{n-1} \Sigma}\nu_{ij}(\hat\bv,\hat\bv_*, u,u_*)
				f_j(\hat\bv_*,u_*)du_*\,d\hat\bv_* {\mbox {  for all  }  {(t,\bx)}  \in [0,+\infty]  \times {\Gamma_{\bx}}},
			\end{aligned}
		\end{equation}
		where the first term on the right-hand side describes the proliferation and the second one describes the destruction. 
		Moreover, the variables involved are
		\begin{itemize}
			\item[--] $\hat\bv',\hat\bv_*\in\mathbb{S}^{n-1}$ and $u',u_*\in\Sigma$ are the velocity direction and the activity of individuals of populations $C_i$ and $C_j$, respectively, before the interaction in the proliferative term; 
			\item[--] $\hat\bv\in\mathbb{S}^{n-1}$ and $u\in \Sigma$ are the velocity direction and the activity, respectively, of either the newborn individuals of population $C_i$ after the proliferative interaction or the individuals of population $C_i$  {before} the destructive interaction.
		\end{itemize}
		Also in this case, $\nu_{ij}(\hat\bv,\hat\bv_*, u,u_*)$ represents the interaction frequency of two  {cell populations} $(C_i, C_j)$ 
		with velocities directed along $(\hat\bv,\hat\bv_*)$ and activities $(u,u_*)$, respectively.
		Similarly for the interaction frequency $\mu_{ij}(\hat\bv',\hat\bv_*, u',u_*)$.
		Moreover, the function $\varphi_{ij}(\hat\bv, u; \hat\bv',\hat\bv_*, u',u_*) $ represents the  fraction of newborn $C_i$ cells with 
		activity $u$ and {velocity directed along} $\hat\bv$ after the interaction between {an individual of} $C_i$ with parameters $(\hat\bv',u')$ 
		and {an individual of} $C_j$ with parameters $(\hat\bv_*,u_*)$. 
		Unlike $\beta_i^H$ satisfying condition (\ref{prob-beta}),
		$\varphi_{ij}$ is not a probability density, since it holds
		\begin{equation}
			\iint\limits_{\mathbb{S}^{n-1} \Sigma}
			\varphi_{ij}(\hat\bv, u; \hat\bv',\hat\bv_*, u',u_*) \,du\,d\hat\bv=\theta_{ij}( \hat\bv',\hat\bv_*, u',u_*),
		\end{equation}
		being, in general, $\theta_{ij}( \hat\bv',\hat\bv_*, u',u_*)\not=1$.
		It represents the total expected number of $C_i$ cells generated through the encounters described above. 
		If $\theta_{ij}( \hat\bv',\hat\bv_*, u',u_*) > 1$ then  {these} interactions lead to a population growth, 
		whereas if $\theta_{ij}( \hat\bv',\hat\bv_*, u',u_*) < 1$ {they} lead to a decay.
		As an example, if both $\theta_{ij}( \hat\bv',\hat\bv_*, u',u_*) > 1$ and $\theta_{ji}(\hat\bv',\hat\bv_*,u',u_*) > 1$,  then we are in a situation of mutualistic synergy between population $i$ and~$j$.
		
		Finally, operators $\mathcal Q^i_{jk}$ in \eqref {eq:ops} take into account the fact that encounters between populations 
		$C_j$ and $C_k$ may also lead to a proliferative event relevant to $C_i$, with $i\neq j,k$. 
		They are defined as
		\begin{equation}
			\begin{aligned}
				\mathcal Q^i_{jk}[f_j,f_k](\hat\bv,u)
				= \iint\limits_{\mathbb S^{n-1}\times\mathbb S^{n-1}}\iint\limits_{\Sigma\times\Sigma}
				& \sigma^i_{jk}( \hat\bv',\hat\bv_*, u',u_*) \psi^i_{jk}(\hat\bv, u; \hat\bv',\hat\bv_*, u',u_*)\\[0mm]
				& \, \times f_j(\hat\bv_*,u_*)f_k(\hat\bv',u')du_*\,du'\,d\hat\bv_*\,d\hat\bv'{\mbox { for all }  {(t,\bx)}  \in [0,+\infty]  \times {\Gamma_{\bx}}}.
			\end{aligned}
			\label{eq:Qi}
		\end{equation}
		{Variables involved are
			\begin{itemize}
				\item[--] $\hat\bv',\hat\bv_*\in\mathbb{S}^{n-1}$ and $u',u_*\in\Sigma$ are the velocity direction and the activity of individuals of populations $C_j$ and $C_k$, respectively, before the interaction; 
				\item[--] $\hat\bv\in\mathbb{S}^{n-1}$ and $u\in \Sigma$ are the velocity direction and the activity, respectively, 
				{of} the newborn individuals of  population $C_i$ after the interaction between the individuals of  population $C_j$ and $C_k$.
		\end{itemize}}
		Again, $\sigma^i_{jk}$ are {interaction} frequencies, 
		and $\psi^i_{jk} (\hat\bv, u; \hat\bv_*,\hat\bv',u_*,u')$ the expected fractions of newborn $C_i$ cells having activity $u$ and 
		velocity directed along $\hat\bv$ after the interaction between {individuals of} $C_j$ and $C_k$.
		The total expected number of new cells $C_i$ is given by
		\begin{equation}
			\iint\limits_{\mathbb{S}^{n-1} \Sigma} \psi^i_{jk}(\hat\bv, u; \hat\bv',\hat\bv_*, u',u_*)\,du\,d\hat\bv=\gamma^i_{jk}( \hat\bv',\hat\bv_*, u',u_*).
		\end{equation}
		
		System \eqref{EqKin} describes the evolution of $N$ cellular populations that,	in addition to interacting among themselves, can diffuse in the host medium 
		and spread across the spatial domain of evolution.
		In the next section, we will investigate a proper asymptotic limit of equations \eqref{EqKin}, 
		leading to a closed system of reaction-diffusion equations for the number densities of cell populations $C_1, \dots, C_N$.
		
		
		\section{Diffusive limit of the kinetic system}
		\label{SecDiffLim}
		
		{We now consider system \eqref{EqKin} and investigate an asymptotic regime that allows us to formally derive systems of 
			reaction-diffusion type from kinetic equations in the diffusive limit.
			A rigorous justification of the hydrodynamic limit (in the spirit of results obtained in \cite{zhigun2022novel} and references therein for various cell scenarios) is scheduled as future work.}
		As in classical diffusive limits, already investigated also in the gas dynamics frame 
		\cite{Anwasia-etal-2017, bisi2006reactive, bisi2022reaction}, the dominant process in the evolution is 
		the one associated with the conservative interactions of cells with the host medium, 
		That is much denser than the populations $C_i$.
		In other words, we take a small parameter $\epsilon$, representing the Knudsen number, 
		and assume that conservative interactions are of order $1/\epsilon$, while the non-conservative ones are 
		of order $\epsilon$ and thus much less frequent. 
		Since we are also interested in the effects of non-conservative dynamics, we have to measure time in the same scale, 
		i.e., we rescale the time setting $t'=\epsilon\,t$. In the sequel, the apex will be omitted for simplicity. 
		
		In this regime, the scaled kinetic system \eqref{EqKin} becomes
		\begin{equation}
			\label{EqKinSca}
			\epsilon \, \frac{\pa f_i}{\pa t} {\,+\,\nabla_{\bx}\cdot \left(c_i\,\hat\bv f_i\right){+ \, \epsilon\,\nabla_{\bv}\cdot\left(f_i\,\mathbf S_i \right)}} 
			=\frac{1}{\epsilon}\,\mathcal G^H_i[f_i,f_H]+\epsilon\,\mathcal H_i[\fvec] , \quad i=1,\ldots,N.
		\end{equation}
		{We underline that the $\epsilon$ in the third term on the left-hand side of \eqref{EqKinSca} comes from 
			the transport term \eqref{Esse}, because of the time rescale.}
		Our present scope is to derive, from equations \eqref{EqKinSca}, 
		a closed system of equations for the macroscopic densities $n_i(t,\bx),\,i =1,\ldots,N$, defined in (\ref{ni}). 
		To this aim, we consider a Hilbert expansion of each distribution function $f_i$ in terms of the 
		scaling parameter $\epsilon$, writing
		\begin{equation}
			\label{ExpDis}
			f_i(t,\bx, \hat\bv,u)=f_i^0(t,\bx, \hat\bv,u)+\epsilon \,f_i^1(t,\bx, \hat\bv,u)+\epsilon^2 \,f_i^2(t,\bx, \hat\bv,u)+O(\epsilon^3).
		\end{equation}
		Without loss of generality, as proven in \cite{bisi2006reactive, bisi2022reaction}, we may suppose that the whole mass density of each population is concentrated on 
		the $\epsilon^0$ term, that is
		\begin{equation}
			\iint\limits_{\mathbb{S}^{n-1} \Sigma} f_i^0(t,\bx, \hat\bv,u)\,du\,d\hat\bv=n_i(t,\bx),\qquad\iint\limits_{\mathbb{S}^{n-1} \Sigma} f_i^k(t,\bx, \hat\bv,u)\,du\,d\hat\bv=0,\mbox {for}  k\geq1.
		\end{equation}

		We note that, from equations (\ref{EqKinSca}) with expansions \eqref{ExpDis},
		one {obtains} that {the conservative operators describing interactions with the host population  play the dominant role, namely}
		$$
		\mathcal G^H_i[f_i^0,f_H] = O(\epsilon) .
		$$ 
		This means that, to the first order of accuracy, the distribution is an equilibrium state of the {linear operator} 
		$\mathcal G^H_i[f_i,f_H]$. Unlike in the classical kinetic theory of gases, where equilibria of Boltzmann operators are Maxwellian distributions, 
		the equilibrium states of {conservative operators} for social {and biological} sciences are, in general,
		not explicit, since interaction rules are based on probability transitions and may also take into account non-deterministic
		or random effects \cite{cordier2005kinetic, della2021sir}.
		
		As already pointed out in~\cite{Bellomo-Belloquid}, the explicit shape of the equilibrium distributions is not needed, but it suffices to know that they exist, provided the results stated in Subsection \ref{ConsDyn}. Furthermore, these results allow us to develop the asymptotic analysis of the scaled equations \eqref{EqKinSca}.
		The first step consists of substituting the Hilbert expansion \eqref{ExpDis}
		of the distribution functions into the scaled kinetic equations \eqref{EqKinSca},
		and equating the same order terms in $\epsilon$.
		We obtain
		{
			\begin{eqnarray}
				& & \mathcal {G}_i^H[f_i^0,f_H] = 0,
				\label{eq:e0}
				\\[2mm]
				& & \nabla_{\bx}\cdot \left(c_i\,\hat\bv f_i^0\right)= \mathcal {G}_i^H[ f_i^{1},f_H],
				\label{eq:e1}
				\\[2mm]
				& &  \frac{\pa f_i^0}{\pa t} \!+\! \nabla_{\bx}\cdot \left(c_i\,\hat\bv f_i^1\right)
				+\nabla_{\bx}\cdot \left(c_i\,\hat\bv f_i^0\right){+\,\nabla_{\bv}\cdot\left(f_i^0\,\mathbf S_i \right)}=  \mathcal {G}_i^H[ f_i^{2},f_H] \!+\! \mathcal H_i[\fvec^0] .
				\label{eq:e2}
		\end{eqnarray}}
		From equation \eqref{eq:e0}, we determine $f_i^0$. 
		Preliminarily, integrating the equation with respect to $u$, and recalling the expression for $\mathcal G_i$ given in \eqref{ConsOp} we may write 
		
		$$
		\begin{aligned}
			\int\limits_{\Sigma}\mathcal G_i^H[f_i^0,f_H](t,\bx,&\hat\bv,u)\,du
			= \int\limits_{\mathbb S^{n-1}}\iiint\limits_{\Sigma^H\times\Sigma\times\Sigma}\  \Big[\eta_{i}^H(u',u_*)\beta_i^H(u; u',u_*)\, f_i^0(\hat\bv',u') \\[0mm]
			&    -\eta_{i}^H(u,u_*)\beta_i^H(u'; u,u_*)\,f_i^0(\hat\bv,u) \Big]\,f_H(u_*)\,du'\,du\,d u_*\,d\hat\bv' =0,
		\end{aligned}
		$$
		that can be written as
		\begin{equation}\label{Fi0}
			\int\limits_{\mathbb S^{n-1}}\ \left( \mathcal F^0_i(\hat\bv')-\mathcal F^0_i(\hat\bv) \right)\,d\hat\bv'=0, \qquad \quad \forall\, \hat{\bv} \in \mathbb S^{n-1}\,,
		\end{equation}
		with
		\begin{equation}
			\mathcal F^0_i(\hat\bv):=\iiint\limits_{\Sigma^H\times\Sigma\times\Sigma}\  \eta_{i}^H(u,u_*)\beta_i^H(u'; u,u_*)\,f_i^0(\hat\bv,u) \,f_H(u_*)\,du\,du'\,du_*.
		\end{equation}
		Then, we may observe that \eqref{Fi0} implies that functions $\mathcal F^0_i(\hat\bv)$ are constant in $\hat\bv$. 
		At this point, keeping in mind  Assumption \ref{Assu} and using the result of Lemma \ref{Lem},
		we infer that {$f_i^0$ is explicitly given as}
		\begin{equation}
			f_i^0(t,\bx,u)= n_i(t,\bx)\,M_i(u).
		\end{equation}
		
		Next, from equation \eqref{eq:e1} we determine $f_i^1$.
		Substituting $f_i^0$ in \eqref{eq:e1} leads to
		\begin{equation}
			\label{eq:e1Sec_}
			{M_i(u)\,\hat\bv\cdot\nabla_{\bx}\left(\,c_i\, n_i\right)} = \mathcal {G}_i^H[ f_i^{1},f_H].
		\end{equation}
		As discussed in the proof of Lemma \ref{Lem}, for the solvability of equation (\ref{eq:e1Sec_}) it is necessary and sufficient that the integral over $\hat\bv$ and $u$ of the term on the left-hand side is null.
		This requirement is trivially fulfilled, provided that the functions $c_i(t,\bx,u)$ are sufficiently regular. 
		With the application in mind (that will be described in Section 5), we suppose that there exists $\widetilde{c}_i(t,\bx)$ such that $$ c_i(t,\bx,u)=u\, \widetilde{c}_i(t,\bx).$$ 
		This assumption is reasonable from the biological point of view since it means that the speed of the cell, namely its movement rate in the considered tissue, is proportional to the cellular activity.
		{This proportionality can be assumed as long as the activity considered in the specific application enhances cell speed rather than reducing it. We emphasize that the model is presented here in a general setting, 
			with the specific activity chosen \emph{ad hoc} in Section \ref{SecApp}, devoted to applications.} 
		The tilde on $\widetilde{c}_i(t,\,\bx)$ will be omitted in the sequel and \eqref{eq:e1Sec_} will be written as
		\begin{equation}
			\label{eq:e1Sec}
			{u\,M_i(u)\,\hat\bv\,\cdot\nabla_{\bx}\left(\,c_i\, n_i\right)} = \mathcal {G}_i^H[ f_i^{1},f_H].
		\end{equation}
		With this choice for the speed, from Lemma \ref{Lem} we can conclude that it is possible to recover $f_i^1$. 
		Indeed, let us consider the unique solution $\mathbf k_i$ of the equation $\mathcal{G}_i^H[\mathbf k_i, f_H]=\hat\bv\, u\,M_i(u)$.
		Then, equation \eqref{eq:e1Sec} becomes, by linearity,
		\begin{equation}
			\label{eq:e1Ter}
			\mathcal{G}_i^H[{\nabla_{\bx}\left(\,c_i\, n_i\right)}\cdot \mathbf k_i-f_i^{1}, f_H]= 0 .
		\end{equation}
		Therefore, the density $ c_i\, {\nabla_{\bx}\, n_i\,}\cdot \mathbf k_i-f_i^{1}$ is an equilibrium for the {linear operator} $\mathcal{G}_i^H$,
		{hence, repeating the same argument as above, $f_i^1$ is explicitly given as
			\begin{equation}
				f_i^{1}={\nabla_{\bx}\left(\,c_i\, n_i\right)}\cdot \mathbf k_i+h_i^1\,M_i \,,
				\label{eq:fi1}
			\end{equation}
			for a certain function $h_i^1$ such that $f_i^1$ has vanishing density.
			
			Now, the last step is to recover $f_i^{2}$ from equation \eqref{eq:e2}.
			First, using expression \eqref{eq:fi1}, we rewrite equation \eqref{eq:e2} as
			\begin{equation}
				\label{eq:e2Sec}
				\frac{\pa n_i}{\pa t} \, M_i \!+\! u\,c_i\,\,\hat\bv\!\cdot\!\nabla_{\bx} \left({\nabla_{\bx}\left(\,c_i\, n_i\right)}\cdot \mathbf k_i+h_i^1\,M_i  \right){+n_i\,M_i\,\nabla_{\bv}\cdot\mathbf S_i }
				- \mathcal H_i[\underline{\bf n}\,\underline{\bf M}]
				=  \mathcal{G}_i^H[ f_i^{2},f_H] ,
			\end{equation}
			being $\underline{\bf n}\,\underline{\bf M}$ the vector $(n_1\,M_1, n_2\, M_2, \dots, n_N\, M_N)$.
			Considering one more time Lemma \ref{Lem}, the approximation $ f_i^{2}$ can be uniquely recovered 
			from equation \eqref{eq:e2Sec} only if the integral of its left-hand side with respect to $\hat\bv$ and $u$ is null.
			
			This leads to 
			a system of evolution equations for the macroscopic densities $n_i$ in the form
			\begin{equation}\label{eq:e2SecMac}
				\frac{\pa n_i}{\pa t} \!=\,\! \,c_i\,\,\nabla_{\bx}\cdot {
					{\left(\widetilde{\mathcal{\bf D}}_i \cdot\nabla_{\bx}\left(\,c_i\, n_i\right)\right)- n\,\frac{\pa c_i}{\pa t}\,n_i}}
				+\!  \iint\limits_{\mathbb S^{n-1}\times\Sigma} \mathcal H_i[\underline{\mathbf{n}}\,\underline{\bf M}]\,du\,d\hat\bv, {\qquad \text{ for all }  (t, x) \in [0,+\infty]\times \Gamma_{\bx},}
		\end{equation}}
		where $\widetilde{\mathcal{\bf D}}_i$ stands for the tensor
		\begin{equation}\label{ExpDi}
			\widetilde{\mathcal{\bf D}}_i = - \iint\limits_{\mathbb{S}^{n-1}\times\Sigma} u\,\hat\bv \otimes \mathbf k_i\left(\hat\bv,\,u\right) du\, d\hat\bv \,.
		\end{equation}
		
		We recall that the function $\mathbf k_i(\hat\bv,\,u)$ is related to the equilibrium distribution $M_i(u)$ by the equation $\mathcal{G}_i^H[\mathbf k_i, f_H]=\hat\bv\,u\, M_i(u)$, therefore it cannot be made explicit in the general case. In Section~\ref{SecApp} we will show a specific application where the equilibrium $M_i(u)$ and, consequently, $\mathbf k_i\left(\hat\bv,\,u\right)$ and $\widetilde{\mathbf{D}}_i$ will be completely explicit.
		
		Anyway, we can prove that diagonal entries of the diffusion matrix $\widetilde{\mathbf{D}}_i$ are always non-negative, as physically expected. Indeed, 
		$$
		\begin{array}{ccl}
			\widetilde{\mathcal{\bf D}}_i &= & \displaystyle -\, \iint\limits_{\mathbb{S}^{n-1}\times\Sigma} u\,\hat\bv \otimes \mathbf k_i\left(\hat\bv,\,u\right) du\,d\hat\bv\,= - \iint\limits_{\mathbb{S}^{n-1}\times\Sigma} u\,\hat\bv \,M_i(u)\otimes \mathbf k_i\left(\hat\bv,\,u\right)  \frac{du\,d\hat\bv}{M_i(u)} \vspace*{0.2 cm}\\
			& = & \displaystyle -\iint\limits_{\mathbb{S}^{n-1}\times\Sigma} \mathcal{G}_i^H[\mathbf k_i, f_H]\left(\hat\bv,\,u\right) \,\otimes \mathbf k_i\left(\hat\bv,\,u\right)  \frac{du\,d\hat\bv}{M_i(u)}\,,
		\end{array}
		$$
		thus its diagonal entries are non--negative recalling (\ref{ineq-sign}).

		
		\section{Kinetic model leading to reaction-diffusion systems with cross-diffusion}
		\label{SecDiffLim2}

		In this section, we investigate more refined interaction dynamics with respect to those considered in Section \ref{SecKin}, which provided
		the kinetic equations \eqref{EqKin} and their diffusive limit equations derived in Section \ref{SecDiffLim}. 
		Specifically, in addition to the spatial dependence of the speed of the cells through the a-priori fixed function $c_i(t,\,\bx,\,u)$, we assume now that the cells may change their orientation 
		in relation to the other cells.    
		
		The fact that the run-and-tumble movement of a cell can be influenced by a bias originating from an external field 
		was originally modeled by means of a turning operator in \cite{alt1980biased}. 
		Then, in the following works \cite{othmer1988models,othmer2000diffusion,othmer2002diffusion}, this external bias was modeled as the gradient of a chemotactic attracting substance.
		
		In our model, we assume that each cell may change its orientation when moving,
		depending on the macroscopic density of the other $N-1$ populations. 
		In other words, each cell adjusts its orientation depending on the
		concentration of other cellular populations around its spatial neighborhood.
		This is the new effect introduced in the kinetic description.
		
		Moreover, in view of considering different orders of dominance for every process involved in the dynamics, we suppose that the re-orientation process is faster than the non-conservative interactions, but slower than the conservative ones.
		With reference to the scaled equations \eqref{EqKinSca},
		a new term will appear describing the orientation of the cells,
		and,  as for the motion term $c_i\hat{\bv}\cdot\nabla_{\bx}f_i $, 	this is assumed of order $\epsilon^0$.
		Thus, we consider the system of dimensionless equations
		\begin{equation}
			\label{EqKinSca2}
			\epsilon \, \frac{\pa f_i}{\pa t}{\,+\,\nabla_{\bx}\cdot \left(c_i\,\hat\bv f_i\right){+\,\epsilon\,\nabla_{\bv}\cdot\left(f_i\,\mathbf S_i \right)}} =\frac{1}{\epsilon}\,\mathcal G^H_i[f_i,f_H]+ \sum_{\substack{j=1 \\ j\ne i}}^N \mathcal L_{ij}[f_i]+\epsilon\,\mathcal H_i[\fvec], \quad i=1,\ldots,N.
		\end{equation}
		The terms $\mathcal{L}_{ij}[f_i]$, with $j\ne i$, represent the turning operators and take the form
		\begin{equation}
			\mathcal{L}_{ij}[f_i](t,\bx,\,\hat\bv,u)=\int_{\mathbb S^{n-1}} T_{ij}( \hat\bv;t,\bx,\hat\bv')f_i(t,\bx,\hat\bv',u)d\hat\bv',
			\label{eq:turn}
		\end{equation}
		with the turning { rates} $T_{ij}$ given by
		\begin{equation}
			{
				T_{ij}( \hat\bv;t,\bx,\hat\bv')=\ \lambda_{ij}(t,\bx,u)\,\widetilde T_{ij}(\hat\bv;t,\bx,\hat\bv'), \quad \widetilde T_{ij}(\hat\bv;t,\bx,\hat\bv')= p_{ij}\hat{\bv}\cdot\hat{\bv}'(\hat{\bv}'\cdot\nabla_{\bx}n_j(t,\bx))\,.}
		\end{equation}
		These  { rates} $T_{ij}$ describe the re-orientation of the cells of population $i$ from $\hat\bv'$ to $\hat\bv$ as depending on the actual
		orientation towards the concentration gradient of the {population $j$} and are
		influenced by a { turning coefficient}   {$\lambda_{ij}>0$} that may depend also on time, space and activity. {The factors $p_{ij}$ {are constant and equal to either} $1$ or $-1$}, depending on the fact that the action of the $j$-th population on the $i$-th one is of attractive or repulsive type, respectively.
		Indeed, if $\hat\bv' \cdot \nabla_{\bx}n_j(t,\bx) >0$, then  $T_{ij}$ takes its maximum value for $\hat\bv = \hat\bv'$ when {$p_{ij}=1$} and for $\hat\bv = -\, \hat\bv'$ when  {$p_{ij}=-1$}; the opposite holds when  $\hat\bv' \cdot \nabla_{\bx}n_j(t,\bx) <0$. This means that the choice  {$p_{ij}=1$} forces the cells to move in a direction close to that of $\nabla_{\bx}n_j(t,\bx)$, while the option  {$p_{ij}=-1$} pushes the cells in the opposite direction.
		
		{ We remark that in usual formulations (see \cite{eckardt2024mathematical,othmer2002diffusion}), 
			the turning operator appears in an expanded form. In our case, interactions at the leading order between the bacterial population and the host environment are included in the operators $\mathcal G_i^H$, and the turning operators $\mathcal L_{ij}$ are smaller perturbations.}
		
		Considering the new scaled equations \eqref{EqKinSca2} as the starting point, we apply the same asymptotic procedure implemented in Section \ref{SecDiffLim} 
		to obtain a reaction-diffusion system for macroscopic densities. 
		Accordingly, let us suppose also that Assumption \ref{Assu} holds, along with the result stated in Lemma \ref{Lem}
		for the conservative operator $\mathcal G^H_i[f_i,f_H]$. 
		
		Then, we insert the expansions \eqref{ExpDis} in the scaled equations \eqref{EqKinSca2} and equal the same order terms, getting
		\begin{eqnarray}
			& & \mathcal {G}_i^H[f_i^0,f_H] = 0,
			\label{eq:e0_L}
			\\[2mm]& & {\nabla_{\bx}\cdot \left(c_i\,\hat\bv f_i^0\right)}  = \mathcal {G}_i^H[ f_i^{1},f_H]+\sum_{\substack{j=1 \\ j\ne i}}^N \mathcal L_{ij}[f_i^0],
			\label{eq:e1_L}
			\\[0mm]
			& &  {\frac{\pa f_i^0}{\pa t} \!+\! {\nabla_{\bx}\cdot \left(c_i\,\hat\bv f_i^1\right)
					+\nabla_{\bx}\cdot \left(c_i\,\hat\bv f_i^0\right){+\,\nabla_{\bv}\cdot\left(f_i^0\,\mathbf S_i \right)}}
				=  \mathcal {G}_i^H[ f_i^{2},f_H]+ \sum_{\substack{j=1 \\ j\ne i}}^N \mathcal L_{ij}[f_i^1] \!+\! \mathcal H_i[\fvec^0]
				.}
			\label{eq:e2_L}
		\end{eqnarray}
		With the same assumptions introduced in the previous section, we can again recover the zeroth-order term of the 
		distribution functions from equations \eqref{eq:e0_L} as
		\begin{equation}
			f_i^0(t,\bx, \hat\bv,u)= n_i(t,\bx)\,M_i(u),
		\end{equation}
		where $M_i(u)$ is the equilibrium distribution for the linear operator $\mathcal {G}_i^H$ supposed to exist in Assumption \ref{Assu}.
		Next, from equation \eqref{eq:e1_L}, we obtain
		\begin{equation}
			\label{eq:e1Sec_L}
			\Bigg({\hat\bv\cdot\nabla_{\bx}\left(\,c_i\, n_i\right)}-\sum_{\substack{j=1 \\ j\ne i}}^N \mathcal L_{ij}[ {n_i}] \Bigg) M_i
			= \mathcal {G}_i^H[ f_i^{1},f_H],
		\end{equation}
		where
		$$
		\mathcal L_{ij}[ {n_i} ] =  n_i(t,\bx)\, {p_{ij}} \lambda_{ij}(t,\bx,u)\,\int_{\mathbb S^{n-1}}  \hat{\bv}\cdot\hat{\bv}'(\hat{\bv}'\cdot\nabla_{\bx}n_j(t,\bx))\,d\hat\bv'\,.
		$$	
		We observe that, also in this case, the integral over $\hat\bv$ and $u$ of the term on the left-hand side of equation \eqref{eq:e1Sec_L} is null,  
		provided that the functions $ c_i(t,\,\bx,u)$ and $\lambda_{ij}(t,\,\bx,\,u)$ are sufficiently regular. 
		Coherently with choices performed in the previous Section, we take $ c_i(t,\,\bx,u)= u\,\widetilde{c}_i(t,\,\bx)$ and 
		$\lambda_{ij}(t,\,\bx,\,u)=u\,\widetilde{\lambda}_{ij}(t,\,\bx)$ (tildes will be omitted in the sequel).
		Due to this, we find it convenient to write equation \eqref{eq:e1Sec_L} as
		\begin{equation}\label{eq:e1Sec_L_2}
			\Bigg({\nabla_{\bx}\left(\,c_i\, n_i\right)}-\sum_{\substack{j=1 \\ j\ne i}}^N \widetilde{\boldsymbol {\mathcal L}}_{ij}[ {n_i}] \Bigg)\cdot\hat\bv\, u\, M_i
			= \mathcal {G}_i^H[ f_i^{1},f_H],
		\end{equation}
		with
		$$
		\widetilde{\boldsymbol {\mathcal L}}_{ij}[ {n_i}]= {p_{ij}} n_i(t,\bx)\, \lambda_{ij}(t,\bx)\,\int\limits_{\mathbb S^{n-1}}  \hat{\bv}'(\hat{\bv}'\cdot\nabla_{\bx}n_j(t,\bx))\,d\hat\bv'\,$$ 
		$$\quad\quad= {p_{ij}} n_i(t,\bx)\, \lambda_{ij}(t,\bx)\,\nabla_{\bx}n_j(t,\bx)\int\limits_{\mathbb S^{n-1}}  (\hat{\bv}')_k^2\,d\hat\bv'\,,
		$$
		denoting by $(\hat{\bv}')_k$ the $k$-th component of $\hat\bv$.
		Thus, we determine the first-order approximation of the distribution functions as
		\begin{equation}
			f_i^{1}=	\Bigg({\nabla_{\bx}\left(\,c_i\, n_i\right)}-\sum_{\substack{j=1 \\ j\ne i}}^N\widetilde{\boldsymbol {\mathcal L}}_{ij}[ {n_i}]\Bigg)
			\cdot \mathbf k_i+h_i^1\,\,M_i,\end{equation}
		still being $\mathbf k_i\left(\hat\bv,\,u\right)$ the unique solution of $\mathcal G_i^H[\mathbf k_i,f_H]=\hat\bv\,u\,M_i(u)$.

		Lastly, let us now consider equation \eqref{eq:e2_L}, that can be rewritten as
		\begin{equation}\label{eq:e2Sec_L}
			\begin{aligned}
				\frac{\pa n_i}{\pa t} M_i & \!+\! u\,c_i\,\,\hat\bv\!\cdot\!\nabla_{\bx} 
				\Bigg(	\Bigg({\nabla_{\bx}\left(\,c_i\, n_i\right)}-\sum_{\substack{j=1 \\ j\ne i}}^N\widetilde{\boldsymbol {\mathcal L}}_{ij}[ {n_i}]\Bigg)
				\cdot \mathbf k_i+h_i^1\,\,M_i\Bigg) {+n_i\,M_i\,\nabla_{\bv}\cdot\mathbf S_i }
				\\[2mm]
				& -\mathcal H_i[\fvec^0] - \sum_{\substack{j=1 \\ j\ne i}}^N \mathcal L_{ij}[f_i^1]
				=  \mathcal {G}_i^H[ f_i^{2},f_H] .
			\end{aligned}
		\end{equation}
		The solvability condition {stated} in Lemma \ref{Lem}, {enabling the derivation of}
		an expression for {the second-order approximation} $f_i^{2}$,
		{requires} the integral over $\hat\bv$ and $u$ of the term on the left-hand side of equation  \eqref{eq:e2Sec_L} to be zero. 
		By imposing such constraint, 
		we derive reaction-diffusion equations for $n_i$, $i=1, \dots, N$, where cross-diffusion terms {emerge}:
		\begin{equation}
			\label{eq:e2SecMac_L}
			\frac{\pa n_i}{\pa t} \!=\,\! \,c_i\,\,\nabla_{\bx}\cdot \left(	{\left(\widetilde{\mathcal{\bf D}}_i \cdot\nabla_{\bx}\left(\,c_i\, n_i\right)\right)- n\,\frac{\pa c_i}{\pa t}\,n_i}-\widetilde{\boldsymbol{\chi}}_i \cdot n_i\,\sum_{\substack{j=1 \\ j\ne i}}^N {p_{ij}} \, \lambda_{ij}\nabla_{\bx} n_j \! \right)
			\!+\!  \iint\limits_{\mathbb S^{n-1}\times\Sigma} \mathcal H_i[\underline{\mathbf{n}}\,\underline{\bf M}]\,du\,d\hat\bv,
		\end{equation}
		{ for all   $(t, x) \in [0,+\infty]\times \Gamma_{\bx}$,} with $\widetilde{\mathbf{D}}_i$ as in \eqref{ExpDi} and \begin{equation}\label{ExpXi}
			\displaystyle\widetilde{\boldsymbol{\chi}}_i=\widetilde{\mathbf{D}}_i\,\int\limits_{\mathbb S^{n-1}}  (\hat{\bv}')_k^2\,d\hat\bv'\,,
		\end{equation}
		while the integral of the last term on the left-hand side of \eqref{eq:e2Sec_L} vanishes due to the shape of the turning operator \eqref{eq:turn}.
		Equations \eqref{eq:e2SecMac_L} are reaction-diffusion equations 
		describing cross-diffusion, as a consequence of introducing the turning operators \eqref{eq:turn}
		in the kinetic equations \eqref{EqKinSca2}, accounting for the re-orientation of the cells. Diagonal entries of matrix $\widetilde{\mathbf{D}}_i$ are non--negative (see the proof at the end of Section \ref{SecDiffLim}); thus the sign of cross-diffusion effects is determined by factors ${p_{ij}}$ appearing in the turning operators that, as already discussed, are positive in the attractive case and negative in the repulsive case.

		
		\section{Application to bacterial communities on a leaf surface}
		\label{SecApp}
		
		In this section, we apply the  approach outlined in the previous sections to study a concrete problem involving 
		bacterial communities living on a leaf surface. 
		These communities interact with each other, and also with the leaf, experiencing different reproductive or destructive processes and competing for resources.
		Additionally, different interactions between bacteria and their environment occur at various spatial scales. 
		Understanding these scales is crucial for thoroughly interpreting microbial colonization patterns. 
		%
		%
		In the following, we adapt the nomenclature from previous sections to the context of the problem considered in this application and, since we are considering a small portion of a leaf, we assume the space domain ${\Gamma_{\bx}} \subset \mathbb{R}^2$ and the velocity direction $\hat{\bv} \in \mathbb{S}^1$.

		
		We consider as biological setting two bacterial strains, denoted by $C_1$ and $C_2$,
		moving on a leaf surface, also known as the phyllosphere, that presents a diverse and intricate environment where microbial inhabitants 
		contend with fluctuating conditions, including varying resource availability, interactions with other microbes, 
		and exposure to environmental stresses like UV radiation, temperature variations, and dryness. 
		Additionally, the leaf's surface exhibits distinct topography and structural elements like stomata, trichomes, and veins, 
		with each one impacting microbial adaptation in different ways. 
		Understanding the interplay among these factors and their influence on microbial communities within the phyllosphere poses a challenging task. 
		Certain factors may exert more localized effects compared to others, adding complexities to the study \cite{esser2015spatial}.
		
		As previously mentioned, the kinetic theory of active particles takes into account how the state of individual cells may change through interactions among different cells. 
		In the context we are considering here, bacterial cells can interact over ``long distances'', by changing the concentration of solutes, such as nutrients in their environments, thereby influencing fluxes of compounds and metabolites diffusing from cell to cell \cite{franklin2007statistical}. 
		Moreover, specific conditions on the leaf surface may induce direct physical interactions among cells \cite{tecon2018cell}.
		
		{The host population $H$ {consists of} cells of the leaf surface.
			Studies indicate that bacteria are more likely to thrive in higher humidity conditions \cite{burkhardt1999measurements}
			and, for this reason, the activity $u\in\Sigma^H$ of the host cells, with $\Sigma^H=[0,1]$, will be {defined} as the {humidity} level right above the surface.  
			While this condition may vary across different parts of the leaf and over the day, 
			we consider a sufficiently small portion of the leaf {where} the distribution $f_H(u)$ is assumed to remain constant in space and time.
			{Specifically}, for the sake of simplicity we {assume a} uniform distribution on~$\Sigma ^H$, {namely} $ f_H(u)=1.$}
		
		{ {Furthermore}, the activity of bacteria $u$ is represented by their {motility}, {which is supposed to belong} to $\Sigma=[-1,1]$, understood as the effective capacity of cells to move across the leaf surface through mechanisms such as flagellar activity, twitching, or surfactant-mediated sliding.  Experimental studies demonstrate that stronger motility mechanisms lead to faster and more extensive displacement. For instance, surfactant-producing Pseudomonas spp. enable co-swarming on leaf surfaces, thereby enhancing spreading speed \cite{kunzler2024hitching}. Similarly, Pseudomonas syringae mutants deficient in surfactant production show markedly reduced surface spreading, confirming that motility mechanisms directly affect displacement rates \cite{holscher2017sliding}. These findings support our assumption that bacterial speed is proportional to motility in the general model framework, an assumption that will be applied in subsequent paragraphs.}
		
		In our model, interactions of cells with the host environment are of a conservative type, 
		leading to either an increase or decrease in the activity, depending on the dryness or moisture of the surface \cite{harshey2003bacterial}.
		We model such dynamics using the following conservative operator, where, compared to its general form (\ref{ConsOp}), the interaction frequencies $\eta_i^H$ and transition probabilities~$\beta_i^H$ are assumed constant with respect to cellular activity {and velocity}, and denoted by $\overline\eta_{i}$ and~$\overline\beta_i$, respectively.  {In this sense,  $\overline\eta_{i}$ and~$\overline\beta_i$ are the probabilities of a cell to change,  respectively, its activity and velocity.  Then,  we have:}
		\begin{equation}
			\mathcal G_i^H[f_i,f_H](\hat\bv,\,u)=\overline\eta_{i}\overline\beta_i\,\int\limits_{\mathbb S^1} 
			\int\limits_{-1}^1\, \left[ f_i(\hat\bv',\,u')-f_i(\hat\bv,\,u)\right]\,du'\,d\hat\bv',\quad i=1,2,
		\end{equation}
		(the dependence of this and the following operators on $t$ and ${\bf x}$ will be omitted); in this particular case, we have, due to the normalization property \eqref{prob-beta}, $\overline\beta_i=1/4\pi$.
		
		Another key factor of microbial interactions with the phyllosphere is the nourishment of bacteria. 
		Indeed, the survival of bacterial cells on a leaf is related to how nutrients like polysaccharides are accessible to microorganisms, 
		on the one hand,
		and how bacteria can modulate the permeability of the leaf, for example, through the
		production of biosurfactants, {on the other hand}.
		See, for example, references \cite{burch2014hygroscopic,krimm2005epiphytic,schreiber2005plant}.
		For this reason, we introduce another population in our model, namely population $L$, 
		constituted by the cells of the leaf surface where nutrients are prevalent and disposable for bacteria, 
		such as close to specialized epidermal outgrowths as trichomes, 
		above veins, and in epidermal cell grooves \cite{brewer1991functional, morris1997methods}.
		For the population $L$, the activity $u\in[-1,1]$ of the cells represents the nourishing capacity, and the distribution function $f_L$ depends on activity, time, and space. 
		The interaction of cells of population $L$ with the populations $C_1$ and $C_2$ leads to the proliferation of the microbial populations {$C_1$ or $C_2$} and the detriment of cells of population $L$. 
		Thus, these interactions are modeled by the following non-conservative operators,
		\begin{eqnarray}
			\mathcal N_{iL}[f_i,f_L](\hat \bv, u) 
			\!\!\!&\!=\!&\!\!\! \overline\mu_{iL} \varphi_{iL}(\hat\bv, u) \!  \! \int\limits_{\mathbb S^1} \int\limits_{-1}^1 \int\limits_{-1}^1 \!\! f_i(\hat\bv_*,u_*)f_L(u')du_* du' d\hat\bv_*, 
			\  i\!=\!1,2, \hspace*{0.5cm}
			\label{eq:NiL}
			\\
			\mathcal N_{Li}[f_L,f_i](u) 
			\!\!\!&\!=\!&\!\!\!  -\overline\nu_{Li} \, f_L(u) \int\limits_{\mathbb S^1}\,\int\limits_{-1}^1
			{f_i} (\hat\bv',u') \, {du' d\hat\bv'}, \ i=1,2.
			\label{eq:NLi}
		\end{eqnarray}
		Again, with respect to the general operator (\ref{coll_ope}), interaction frequencies are supposed constant ($\overline \mu_{iL}$ and $\overline \nu_{Li}$), and the expected fraction of newborn cells $C_i$ (represented by the function $\varphi_{iL}$) is assumed independent of the pre--interaction parameters.
		
		Observations in laboratory \cite{esser2015spatial} {indicate} that {the frequent} co-aggregation of cells 
		from different strains suggests that two populations may somehow either facilitate each other or exploit resources similarly in the phyllosphere.
		We incorporate these dynamics in our model by introducing, in the kinetic equation for population $L$, a term that describes how such interactions may increase the number of available nourishing cells, namely
		\begin{equation}
			\mathcal Q_{12}^L[f_1,f_2](\hat \bv, u) \!=\! \overline\sigma_{12}^L \,\psi_{12}^L(\hat\bv, u) \!\!
			\iint\limits_{\mathbb S^1\times\mathbb S^1}\int\limits_{-1}^1 \! \int\limits_{-1}^1 \!\!
			f_1(\hat\bv_*,u_*) f_2(\hat\bv',u') \, du_* du' d\hat\bv_* d\hat\bv', \ \ i=1,2,
		\end{equation}
		where  {the interaction frequency $\overline{\sigma}_{12}^L$ is constant and} $\psi^L_{12} (\hat\bv, u)$, as already explained in Section \ref{SecKin}, denotes the expected fraction of newborn $L$ cells having activity $u$ and velocity direction $\hat\bv$ after the interaction between a $C_1$ and a $C_2$ cell and is assumed independent of pre--interaction parameters.

		Regarding competition among individuals, two different dynamics will be considered, namely
		interference and exploitation. 
		Interference occurs when one species actively excludes others, 
		often through mechanisms like antibiosis, where toxic compounds are produced. 
		{An example of this is} described in \cite{pusey2011antibiosis}. 
		Exploitation, on the other hand, arises from the competition for shared resources, like nutrients or space. In the present scenario, competition may result in compromising population growth due to resource limitations
		\cite{remus2012variation}. 
		In our model we suppose that {$C_2$-}population is more aggressive than the other, 
		and we {thus} define the {following destructive} operators {associated to interference and exploitation} as   
		\begin{eqnarray}
			\mathcal N_{12}[f_1,f_2](\hat \bv, u) &\!\!=\!\!& -\overline\nu_{12}\,f_1(\hat\bv,u)\int\limits_{\mathbb S^1}\,\int\limits_{-1}^1f_2(\hat\bv',u')\,du'\,d\hat\bv',
			\\
			\mathcal N_{22}[f_2,f_2](\hat \bv, u) &\!\!=\!\!& - \overline\nu_{22}\,f_2(\hat\bv,u)\int\limits_{\mathbb S^1}\,\int\limits_{-1}^1f_2(\hat\bv',u')\,du'\,d\hat\bv',
		\end{eqnarray}
		{where the operator $\mathcal N_{12}[f_1,f_2]$ describes the interspecific competition resulting in the destructive effect of the $C_2$-population on the $C_1$-population, 
			whereas the operator $\mathcal N_{22}[f_2,f_2]$ describes the intraspecific competition within the $C_2$-population.  {Constant interaction rates are denoted by $\overline{\nu}_{12}$ and $\overline{\nu}_{22}$, respectively.}} 
		
		Moreover, the natural death of only the aggressive strain will be considered in the model. Thus we have
		\begin{equation}
			\mathcal J_{2}[f_2](\hat\bv,u)=- \overline \tau_{2}\,f_2(\hat\bv,u),
		\end{equation}
		with $\overline \tau_{2}$ positive constant.
		For what concerns the movement of bacteria, as already outlined, it is reasonable to assume that the speed
		is proportional to the motility; for this reason, we adopt the modeling introduced in the previous sections. 
		
		Furthermore, we suppose that the population speed depends on space and time through the macroscopic densities of the two strains $C_1$ and $C_2$,
		{so that the cellular velocity is given as} 
		$$
		\bv = \hat\bv\,u\,c_i\left(n_1(t,\,\bx),n_2(t,\,\bx)\right), \quad i=1,2 .
		$$ 
		{ In addition, we suppose that, in the present setting, the derivatives of functions $c_i$ with respect to the macroscopic densities are of order $\epsilon$. More precisely, we introduce vector functions $\mathbf{s}_i$ of order $O(1)$ such that 
			$$
			\nabla_{(n_1,n_2)} \, c_i(n_1,n_2) = \epsilon \, \mathbf{s}_i(n_1,n_2).
			$$
			Under this scaling, the terms involving spatial gradients and temporal derivatives of $c_i$
			are of higher order and can therefore be neglected in the derivation of the macroscopic equations.
			Including these terms would, of course, provide a more rigorous and complete analysis of the model, 
			but we leave this aspect open for future developments.}
		
		As already mentioned, bacteria can produce signaling by changing the concentration
		of solutes in the environment, {up to}
		a distance of approximately ten times their cell’s diameter. We suppose that this type of interaction may induce re-orientation of the cells \cite{othmer2000diffusion}, depending on the orientation of the spatial gradient of the number density of the other population {restricted to the short range, in order to disregard spatially nonlocal terms}. Thus, we consider the following turning operators
		\begin{eqnarray}
			\mathcal{L}_{ij}[f_i](t,\bx,\,\hat\bv,u) 
			&\!\!=\!\!&\int\limits_{\mathbb{S}^1} \lambda_{ij} \,{p_{ij}} \widehat{\bv}
			\cdot\hat{\bv}'(\hat{\bv}'\cdot\nabla_{\bx}n_j(t,\bx)) \, f_i(t,\bx,\hat\bv',u)d\hat\bv', 
			\label{eq:opLij} \\[2mm]
			&  &  \mbox{for} \quad (i,j)\in\{(1,2),(2,1)\}. 
			\nonumber
		\end{eqnarray}
		Inspired by models proposed in the literature \cite{othmer2002diffusion}, we assume that the  { turning coefficients}  are, indeed, 
		functions of the macroscopic densities of the involved populations. Also in this case, the assumption that the turning movement of bacteria depends proportionally on the motility holds. Consequently, we have that the  { turning coefficient}  is}
	$$ {u\,\lambda_{ij}\left(n_1(t,\,\bx),n_2(t,\,\bx)\right), \quad (i,j)\in\{(1,2),(2,1)\}.}
	$$
	Research on the evolution of bacterial populations on leaf surfaces suggests that interactions among bacteria often occur at small spatial scales, in contrast to interactions between bacteria and the environment. 
	Consequently, it is reasonable to consider interactions of the two populations with the host environment as the dominant processes, 
	while interactions among bacteria as slower processes.
	
	With reference to the asymptotics developed in the previous sections, and to the small parameter $\epsilon$ representing a ratio between a microscopic and a macroscopic scale,	we assume that interactions of the two populations with the host environment are of order $1/\epsilon$, whereas interactions among bacteria are of order $\epsilon$.
	We also suppose that the cooperative interaction between  {the two bacteria populations,  described by the term $\mathcal Q_{12}^L[f_1,f_2]$, and the turning operators $\mathcal L_{ij}[f_i,f_j]$ are} of order~$\epsilon^0$.
	
	Thus, the kinetic system for distribution functions $f_1,\,f_2,\,f_L$ describing the biological setting under investigation is the following:
	\begin{equation}
		\begin{aligned}\label{SistKinBio2}
			\epsilon\,\frac{\pa f_1}{\pa t} \!+\! u\,c_1\hat\bv \!\cdot\! \nabla_{\bx}f_1+&{ \epsilon \, u\,f_1\,\mathbf s_1\cdot\left(\hat\bv\left[\nabla_{\bx} n_1,\nabla_{\bx}\,n_2\right]\right) +\epsilon^2\nabla_{\bv}\cdot\left(f_1\,u\,\hat\bv\,\mathbf s_1\cdot\left(\frac{\pa n_1}{\pa t},\frac{\pa n_2}{\pa t}\right)\right)}\\[2mm]\quad
			\! &= \frac{1}{\epsilon}\,\mathcal G^H_1[f_1,f_H]+\mathcal L_{12}[f_1, f_2]
			\!+\! \epsilon \Big( \mathcal N_{1L}[f_1, f_L] \!+\! \mathcal N_{12}[f_1,f_2]\Big)
			\\[2mm]
			\epsilon\,\frac{\pa f_2}{\pa t} \!+\!  u\,c_2\hat\bv \!\cdot\! \nabla_{\bx}f_2+&{ \epsilon \, u\,f_2\,\mathbf s_2\cdot\left(\hat\bv\left[\nabla_{\bx} n_1,\nabla_{\bx}\,n_2\right]\right) +\epsilon^2\nabla_{\bv}\cdot\left(f_2\,u\,\hat\bv\,\mathbf s_2\cdot\left(\frac{\pa n_1}{\pa t},\frac{\pa n_2}{\pa t}\right)\right)}\\[2mm]\quad
			&= \frac{1}{\epsilon}\,\mathcal G^H_2[f_2,f_H]+\mathcal L_{21}[f_2, f_1]
			\!+\! \epsilon \Big(\mathcal N_{2L}[f_2, f_L] +\, \mathcal N_{22}[f_2,f_2]\!+\!  \mathcal J_{2}[f_2] \Big)
			\\[2mm]
			\epsilon\,\frac{\pa  f_L}{\pa t}
			&= \epsilon \Big (\mathcal N_{L1}[ f_L,f_1] + \mathcal N_{L2}[ f_L,f_2] \Big) + \mathcal Q_{12}^L[f_1,f_2].
		\end{aligned}
	\end{equation}
	
	In addition, we assume that during interactions between bacteria and nourishing cells, the rate of consumption of nourishing cells is much higher than the rate of proliferation. 
	In other words, for the interaction rates appearing in operators (\ref{eq:NiL}) and (\ref{eq:NLi}), we assume that
	\begin{equation}\label{relPars}
		\overline\nu_{Li}=\frac{1}{\epsilon}\,\overline \mu_{iL}\theta_{iL},\quad i=1,2,
	\end{equation}
	recalling that 
	$$
	\theta_{iL} =\int\limits_{\mathbb{S}^1}\int\limits_{-1}^1 \varphi_{iL}(\hat \bv, u)\,du\,d\hat\bv,\quad i=1,2.
	$$
	
	We observe that {operators $\mathcal G_i^H[f_i,f_H]$ satisfy the conditions} required by Lemma \ref{Lem} 
	{and corresponding Assumption \ref{Assu}, if we take} {$M_i(u) = \dfrac{1}{4\, \pi}$.}
	
	Therefore, we can apply the asymptotic procedure described in the previous section to the first two equations in \eqref{SistKinBio2}, 
	obtaining macroscopic equations for the densities of bacterial populations. 
	{In particular, we recover  \begin{equation}
			f_i^0(t,\bx, \hat\bv,u)= n_i(t,\bx)\,M_i(u),\quad i=1,2,
		\end{equation}
		
		\begin{equation}
			f_i^{1}=	\Bigg({\,c_i\,\nabla_{\bx} n_i}-\sum_{\substack{j=1 \\ j\ne i}}^N\widetilde{\boldsymbol {\mathcal L}}_{ij}[ {n_i}]\Bigg)
			\cdot \mathbf k_i+h_i^1\,\,M_i,\quad i=1,2,\end{equation}
		with $$
		\mathbf k_i(\hat\bv,u) = -\frac{\hat\bv\,u}{\overline \eta_i\,4\, \pi,}
		$$ 
		{the unique solution to the equation}
		$\ \mathcal{G}_i^H[\mathbf k_i(\hat\bv,u), f_H]=\hat\bv\,u\,M_i(u)$.

		Lastly, we have
		\begin{equation}\label{eq:e2Sec_L_m}
			\begin{aligned}
				\frac{\pa n_i}{\pa t} M_i & \!+\! u\,c_i\,\,\hat\bv\!\cdot\!\nabla_{\bx} 
				\Bigg(	\Bigg({\,c_i\,\nabla_{\bx} n_i}-\sum_{\substack{j=1 \\ j\ne i}}^N\widetilde{\boldsymbol {\mathcal L}}_{ij}[ {n_i}]\Bigg)
				\cdot \mathbf k_i+h_i^1\,\,M_i\Bigg) {+n_i\,u\,M_i\,\mathbf s_i\cdot\left(\hat\bv\left[\nabla_{\bx} n_1,\nabla_{\bx}\,n_2\right]\right) }
				\\[2mm]
				& -\mathcal H_i[\fvec^0] - \sum_{\substack{j=1 \\ j\ne i}}^N \mathcal L_{ij}[f_i^1]
				=  \mathcal {G}_i^H[ f_i^{2},f_H] ,\quad i=1,2.
			\end{aligned}
	\end{equation}}{By applying the solvability condition of Lemma \ref{Lem} to \eqref{eq:e2Sec_L_m} and computing} the self-diffusion and cross-diffusion tensors as in \eqref{ExpDi} and \eqref{ExpXi}, respectively, we get the following system of macroscopic equations
	{
		\begin{equation}
			\label{SistMac1}
			\begin{aligned}
				\frac{\pa n_1}{\pa t} \!&=\! \,c_1\,\nabla_{\bx}\!\cdot\! \left(c_1\,\widetilde{\mathcal {D}}_1 \nabla_{\bx}\, n_1
				-\widetilde{{\chi}}_1 \, {p_{12}} \, n_1\, \lambda_{12}\nabla_{\bx} n_2 \! \right)+\overline\mu_{1L}\,\theta_{1L}n_1\, n_L-\overline\nu_{12}n_1\,n_2, \\[2mm]
				\frac{\pa n_2}{\pa t} \!&=\!  \,c_2\,\nabla_{\bx}\!\cdot\! \left(c_2\,\widetilde{\mathcal{D}}_2 \nabla_{\bx}\, n_
				-\widetilde{{\chi}}_2\, {p_{21}} \,n_2\,\lambda_{21}\nabla_{\bx} n_1 \! \right)+\overline\mu_{2L}\,\theta_{2L}n_2\, n_L-\overline\nu_{22}\,n_2^2-\overline\tau_{2}\,n_2,
			\end{aligned}
		\end{equation}
		with
		\begin{equation}
			{\widetilde{\mathcal{D}}_i=\dfrac{1}{6\,\overline\eta_i},\quad \widetilde\chi_i=\dfrac{\pi}{6\,\overline \eta_i}.}
			\label{eq:DChi}
	\end{equation}}
	From the third equation of system \eqref{SistKinBio2}, instead, along with relation \eqref{relPars}, we can write
	\begin{equation}
		\mathcal N_{L1}[f_L^0 ,n_1\,M_1]+\,\mathcal N_{L2}[f_L^0 ,n_2\,M_2]+\mathcal Q_{12}^L[ n_1\,M_1 ,n_2\,M_2]=0,
	\end{equation}
	which, integrated over variables $u$ and $\hat\bv$, provides the relation
	\begin{equation}\label{ExpnL}
		n_L=\frac{\overline \sigma_{12}^L\,\gamma_{12}^L\,n_1\,n_2}{\overline \mu_{1L}\theta_{1L}n_1+\overline \mu_{2L}\theta_{2L}n_2} \, ,
	\end{equation}
	with 
	$$
	\gamma^L_{12}=	\int\limits_{\mathbb{S}^{1}} \int\limits_{-1}^1 \psi^L_{12}(\hat\bv, u)\,du\,d\hat\bv.
	$$
	Expression \eqref{ExpnL} can be plugged into equations \eqref{SistMac1}, leading to the following system of reaction-diffusion equations,
	\begin{equation}
		\label{SistMac2.0}
		\begin{aligned}
			\frac{\pa n_1}{\pa t} \!&=\! \, \,c_1\,\widetilde{\mathcal{D}}_1\,\nabla_{\bx}\cdot \left(c_1\,\,\nabla_{\bx}\, n_1\,
			-\pi\,{p_{12}} \,n_1\, \lambda_{12}\,\nabla_{\bx} n_2\right)+\overline \sigma_{12}^L\,\gamma_{12}^L\,\frac{n_1^2\,n_2}{n_1+\beta\,n_2}
			-\overline\nu_{12}n_1\,n_2, \\[2mm]
			\frac{\pa n_2}{\pa t} \!&=\!  \,c_2\,\widetilde{\mathcal{D}}_2\,\nabla_{\bx}\cdot \left(c_2\,\,\nabla_{\bx}\, n_2-\pi\, {p_{21}} \,n_2
			\,  \lambda_{21}\,\nabla_{\bx} n_1\right)+\overline     
			\sigma_{12}^L\,\gamma_{12}^L\,\frac{\beta\, n_1\,n_2^2}{n_1+\beta\,n_2}
			-\overline\tau_2\,n_2-\overline\nu_{22}\,n_2^2,
		\end{aligned}
	\end{equation}
	with the coefficient $\beta$ being given by $\beta=\dfrac{\overline \mu_{2L}\theta_{2L}}{\overline \mu_{1L}\theta_{1L}}$, having collected $\widetilde{\mathcal D}_i$ in virtue of \eqref{eq:DChi}.

	At this point, we consider it convenient to perform the further time scaling leading to a normalization of the equations.
	Accordingly, we set $\widetilde t=\overline\nu_{12}\,t$. 
	Thus, redefining the coefficients
	\begin{equation}
		\begin{aligned}
			\mathcal D_i = \frac{\widetilde{\mathcal{D}}_i}{\overline\nu_{12}},  \quad
			\lambda_i = \pi\,{p_{ij}} \,\lambda_{ij}\quad \ \ \mbox{for} \ \ i,j=1,2,\,i\neq j
			\\[2mm]
			\zeta=\frac{\overline \sigma_{12}^L\,\gamma_{12}^L}{\overline\nu_{12}},\quad
			\tau=\frac{\tau_2}{\overline\nu_{12}},\quad
			\nu=\frac{\overline\nu_{22}}{\overline\nu_{12}} \, .
		\end{aligned}
		\label{eq:DDD}
	\end{equation}
	Consequently, the reaction-diffusion system \eqref{SistMac2.0} rewrites as
	\begin{equation}
		\label{SistMac2.1}
		\begin{aligned}
			\frac{\pa n_1}{\pa t} \!&=\! \, \,c_1\,{\mathcal{D}}_1\,\nabla_{\bx}\cdot \left(c_1\,\,\nabla_{\bx}\, n_1\,
			-\lambda_{1}\,n_1\,\nabla_{\bx} n_2\right)+\frac{\zeta\,n_1^2\,n_2}{n_1+\beta\,n_2} -n_1\,n_2 \, , 
			\\[2mm]
			\frac{\pa n_2}{\pa t} \!&=\!  \,c_2\,{\mathcal{D}}_2\,\nabla_{\bx}\cdot \left(c_2\,\nabla_{\bx}\, n_2\,
			-\lambda_{2}\,n_2\,\nabla_{\bx} n_1\right)+\frac{\zeta\,\beta\, n_1\,n_2^2}{n_1+\beta\,n_2}
			-\tau\,n_2-\nu\,n_2^2 \, .
		\end{aligned}
	\end{equation}
	{We couple the above system to the following boundary conditions
		\begin{equation}
			\label{SistMac2.1_bc}
			\begin{aligned}
				\left(c_1\,\,\nabla_{\bx}\, n_1\,
				-\lambda_{1}\,n_1\,\nabla_{\bx} n_2\right) \cdot {\pmb \nu}=0 \, , 
				\\[2mm]
				\left(c_2\,\nabla_{\bx}\, n_2\,
				-\lambda_{2}\,n_2\,\nabla_{\bx} n_1\right)\cdot {\pmb \nu} = 0 \, , 
			\end{aligned}
		\end{equation}
		for all  $\pmb x\in \partial{\Gamma_{\bx}}$ and $ {\pmb \nu}$ being the outward normal vector to ${\Gamma_{\bx}}$ at $\partial{\Gamma_{\bx}}$.  Such boundary conditions describe a net zero-flux of each bacterial population at the boundary of the spatial domain.
		
		Ecological interactions and the availability of resources might lead populations to organize and evolve into bacterial clusters.
		Certain Turing patterns constitute a possible approach to describe these kinds of dynamics, and a clear example is represented by the Turing spots.
		Therefore, in the next subsections, we will investigate the conditions under which system \eqref{SistMac2.1},  with boundary conditions \eqref{SistMac2.1_bc}, can develop Turing instability,  eventually leading to bacterial patterns.} 
	
	\subsection{Turing instability: self-diffusion case}
	\label{SecTuring1}
	
	We consider the {simpler} case in which there is no chemotactic motion
	{and therefore only self-diffusion is present in the description.}
	Thus system \eqref{SistMac2.1} reads as
	\begin{equation}
		\label{SistMac2}
		\begin{aligned}
			\frac{\pa n_1}{\pa t} \!&=\! \, \,c_1\,\,\nabla_{\bx}\cdot \left(c_1\,\,{\mathcal{D}}_1\nabla_{\bx}\, n_1\right)+\frac{\zeta\,n_1^2\,n_2}{n_1+\beta\,n_2}
			-n_1\,n_2, \\[2mm]
			\frac{\pa n_2}{\pa t} \!&=\!  \,c_2\,\,\nabla_{\bx}\cdot \left(c_2\,\,{\mathcal{D}}_2\nabla_{\bx}\, n_2\right)+\frac{\zeta\,\beta\, n_1\,n_2^2}{n_1+\beta\,n_2}
			-\tau\,n_2-\nu\,n_2^2 .
		\end{aligned}
	\end{equation}
	To analyze the Turing instability in this scenario, we start by identifying a homogeneous steady state solution, 
	which is the following
	\begin{equation}
		\label{Equi}
		( \overline n_1, \overline n_2)=\left(\frac{\beta \,\tau\,} {(\zeta -1 ) (\beta  \, -\nu )},\frac{\tau}{\beta\, -\nu} \right).
	\end{equation}
	We note that, for the equilibrium state to be biologically meaningful, its components must be positive, 
	and this corresponds to imposing the relation
	\begin{equation}
		\label{CondEqPos}
		\zeta >1,\quad \beta >\nu.
	\end{equation}
	Then, we identify conditions on parameters ensuring the stability of the equilibrium state \eqref{Equi} {with respect to spatially uniform perturbations}. 
	To do this, we first linearize system \eqref{SistMac2} around the equilibrium, 
	writing
	\begin{equation}
		\dfrac{d\mathbf{w}}{dt}\,=\,\mathbf{J}\mathbf{w}\,,
		\label{eq:linW}
	\end{equation}
	with $\mathbf{w}=(n_1-\overline n_1, n_2-\overline n_2)^T$ {being the vector of the deviations with respect to the equilibrium}
	and $\mathbf{J}$ the Jacobian matrix evaluated at the equilibrium
	\begin{equation}
		\mathbf{J}\,=\,\left(\begin{array}{cc}
			\dfrac{\tau \, (\zeta -1)}{\zeta  (\beta-\nu)}& \dfrac{\beta  \,\tau \,}{\zeta  (\nu-\beta \, )}\\
			\\
			\dfrac{\tau (\zeta -1)^2}{\zeta  (\beta \, -\nu)}& \dfrac{\tau(\beta   -\zeta \, \nu)}{\zeta \,(\beta- \nu)}	\end{array}\right).
		\label{eq:Jacob}
	\end{equation}
	Trace and determinant of $\mathbf J$ are, respectively, given by
	\begin{equation}
		\text{Tr}(\mathbf{J})=\frac{\tau \left(\beta -1 +\zeta  (1-\,\nu)\right)}{\zeta  (\beta -\,\nu)}
		\,,\qquad\text{Det}(\mathbf{J})=\frac{{\tau}^2 \, (\zeta -1)}{\zeta  (\beta  -\,\nu)}\, .
	\end{equation}
	We can immediately notice that, holding the existence condition \eqref{CondEqPos}, we have 
	$$
	\text{Det}(\mathbf{J}) > 0 .
	$$
	On the other hand, the homogeneous steady state is stable if also $$
	\text{Tr}(\mathbf{J}) < 0 ,
	$$ 
	that corresponds to imposing
	\begin{equation}\label{zetabar}
		\zeta  > \overline\zeta  \quad \mbox{with} \quad \overline\zeta=\frac{1-\beta} {1-\nu},
		\qquad \mbox{and} \quad \nu>1.
	\end{equation}
	In addition, it can be observed that when parameters satisfy the equality $\zeta=\overline\zeta$,
	the Jacobian matrix \eqref{eq:Jacob} has a pair of purely imaginary eigenvalues,
	and a Hopf bifurcation occurs, see \cite{rionero2019hopf}.
	{On the other hand, under {conditions \eqref{CondEqPos} and \eqref{zetabar}}, the Jacobian matrix has sign pattern
		\[
		{\text{sign}\big( \mathbf{J} \big)} = 
		\begin{bmatrix}
			+ & - \\
			+ & -
		\end{bmatrix}.
		\]
		Therefore, $n_1$ acts as an activator for $n_2$, which is also self-activating, while $n_2$ acts as an inhibitor, and inhibits $n_1$ and itself.
		For our ecological problem, this might describe a parasitism relationship, which is a form of exploitation, 
		where the first bacterial strain is beneficial to the second one, while the latter, on the other hand, harms the former.}

	Now, by applying spatially heterogeneous perturbations to the uniform steady state, i.e. $\mathbf{w} = \mathbf{w}(t,\mathbf{x})$, instead of  \eqref{eq:linW}, we analyse the following linearization of \eqref{SistMac2}
	\begin{equation}
		\dfrac{\partial\mathbf{w}}{\partial t}\,=\,\mathbf{D}\Delta_\mathbf{x}\mathbf{w}+\mathbf{J}\mathbf{w} ,
		\quad \mbox{for} \ (t , \mathbf{x}) \ \ \text{on} \ \ (0,+\infty)\times{\Gamma_{\bx}}\,,
		\label{eq:lindiffW}
	\end{equation}
	considering that no-flux conditions occur at the boundary, that is
	${\boldsymbol{\nu}}\cdot\nabla_{\bf x} \mathbf{w}=0$.
	The diffusion matrix $\mathbf{D}$ in equation \eqref{eq:lindiffW} is given by
	\begin{equation}
		\mathbf{D}\,=\,\left(\begin{array}{cc}
			\widehat{\mathcal D}_{1}&0 \\
			\\
			0 & 	\widehat{\mathcal D}_{2}
		\end{array}\right),
		\label{eq:detD}
	\end{equation}
	where the self-diffusion coefficients are given by
	\begin{equation}
		\widehat{\mathcal D}_{i}={[ c_i(\overline n_1,\,\overline n_2)]^2\,\mathcal D_i}.
	\end{equation}
	{which result from using the expansions 
		$
		c_i(n_1,\, n_2) = c_i + {\nabla_{\bx}} \, c_i \cdot \mathbf{w} + O(\mathbf{\|w\|}^2)
		$
		in system \eqref{SistMac2} and retaining only the linear terms.} 
	{
		We stress that the linearized system \eqref{eq:lindiffW} retains only the leading linear terms and therefore identifies the parameter regions and spatial modes at which a Turing instability first sets in.
		The validity of this calculation is limited to infinitesimal perturbations, namely, if perturbations are small but finite, nonlinear terms become relevant and determine the eventual amplitude, form, and stability of the emergent patterns. A weakly nonlinear analysis, as performed in similar models derived from the kinetic level \cite{bisi2025derivation,BMT}, would give a more accurate description of pattern formation, including different shapes or stability, beyond onset. Since the present work focuses on the kinetic derivation and on establishing a biologically interpretable macroscopic framework, we defer a thorough weakly nonlinear and bifurcation analysis to future work.}
	
	To solve system \eqref{eq:lindiffW}, we use the separation of variables and consider a normal mode expansion in the Fourier series of the unknown function, namely
	\begin{equation}
		\mathbf{w}(\mathbf{x},t)\,=\,\sum_{k\in\mathbb{N}}\xi_ke^{{\ell_k t}} \, \overline{\mathbf{w}}_k(\mathbf{x}) ,
		\label{solutions}
	\end{equation}
	where the eigenfunctions $\overline{\mathbf{w}}_k$ represent independent perturbation modes
	and eigenvalues {$\ell_k$} represent the corresponding linear growth rates.
	Therefore the eigenfunctions $\overline{\mathbf{w}}_k$ solve the time-independent problem
	\begin{equation}
		\begin{cases}
			\Delta \overline{\mathbf{w}}_k+k^2\,\overline{\mathbf{w}}_k\,=\,\boldsymbol{0} , \,&\text {in} \,{\Gamma_{\bx}} ,\\[2mm]
			\boldsymbol{\nu}\cdot\nabla_\mathbf{x}\overline{\mathbf{w}}_k\,=\,0, \,&\text {at} \partial{\Gamma_{\bx}} ,
		\end{cases}
	\end{equation}
	and the scalar coefficients {$\ell_k$} are eigenvalues of the matrix $\mathbf{J}-k^2\mathbf{D}$.
	
	Turing instability is reached whenever the steady state \eqref{Equi} is linearly unstable to spatial perturbations and, consequently, 
	there must exist at least a wavenumber $k$ such that the real part of the corresponding eigenvalue {$\ell_k$} is positive. 
	{The} coefficients {$\ell_k$} are solutions of the dispersion relation
	\begin{equation}
		{\ell}^2+a(k^2){\ell}+b(k^2)\,=\,0,
		\label{autovalori}
	\end{equation}
	with
	\begin{equation}
		a(k^2)\,=\,k^2\text{Tr}(\mathbf{D})-\text{Tr}(\mathbf{J})
	\end{equation}
	and
	\begin{equation} 
		\label{Expb}
		b(k^2)\,=\,\det(\mathbf{D})k^4+g k^2+\det(\mathbf{J})\,,
	\end{equation}
	where the {term} $g$ is a function of both diffusion and reaction coefficients of the problem, 
	that reads
	\begin{equation}
		g=\frac{\tau\,\left(\widehat{\mathcal D}_{1}\,( \zeta\, \nu-\beta)+\widehat{\mathcal D}_{2}(1-\zeta)\right)}{\zeta\,  (\beta-\nu)}.
	\end{equation}	
	To have roots of the dispersion relation (\ref{autovalori}) with positive real part, since from stability of equilibrium it holds $a(k^2)>0$, it must be
	$b(k^2)<0$ for some nonzero $k$. 
	This is satisfied if $b(k^2)$ evaluated in its minimum given by
	\begin{equation}
		k_c^2\,:=-\dfrac{g}{2\det(\mathbf{D})}
		\label{cond1}
	\end{equation}
	is negative, and this is equivalent to imposing
	\begin{equation}
		{4\det(\mathbf{D})\det(\mathbf{J})-g^2}\,<\,0\,.
		\label{cond2}
	\end{equation}
	We observe that expression (\ref{cond1}) {for $k_c^2$} requires that  $g$ must be negative,
	{since $\det(\mathbf{D})>0$, see equations \eqref{eq:DChi}, \eqref{eq:DDD} and \eqref{eq:detD}.}

	The constraints stated above can be expressed by the following conditions,
	\begin{equation}\label{cond2.1}
		\delta\,\,(1-\zeta)+\zeta\,\nu-\beta<0 ,
	\end{equation}
	\begin{equation}\label{cond2.2}
		\delta^2\, (\zeta -1)^2-2 \delta\, (\zeta -1) (\beta\,\,(2 \zeta -1)-\zeta \,\nu)+	\left(\beta-\zeta \,\nu\right)^2>0,
	\end{equation}
	where $\delta=\dfrac{\widehat{\mathcal D}_{2}}{\widehat{\mathcal D}_{1}}$, which can be put together, obtaining
	\begin{equation} 
		\delta>\frac{2 \beta\,  \zeta\,  -\beta\,  ^2-\zeta  \,\nu}{\zeta -1}+2 \sqrt{\frac{\beta\,  \zeta\,( \beta-  \nu)}{\zeta -1}}.
		\label{TurCond}
	\end{equation}
	The Turing condition in our analysis, represented by inequality \eqref{TurCond}, provides the criterion that the system parameters must satisfy to allow the uniform stable equilibrium \eqref{Equi} to become unstable under spatial perturbations. 
	This may lead to the occurrence of a Turing bifurcation. {In particular, from the sign pattern of the Jacobian matrix, we expect the two populations to be in phase, meaning that they will both concentrate over the same spatial locations, see e.g. \cite[Chapter 2]{murray2003mathematical}.
		As we show with some numerical examples in Section \ref{SecSim}, the solutions develop spots, hence describing a tendency of the two bacterial strains to concentrate and co-exist in spatially segregated colonies.
	}

	
	\subsection{Turing instability: {cross-}diffusion case}
	\label{SecTuring2}
	
	We consider now the situation with chemotactic motion leading to cross-diffusion.
	The dynamics are described by the complete system \eqref{SistMac2.1}.
	
	The study of the stability for the {spatial} homogeneous equilibrium $(\overline n_1,\overline n_2)$ is analogous to the one 
	{developed} in Subsection \ref{SecTuring1}. 
	As for the Turing instability, instead, the diffusion matrix results in being
	\begin{equation}\label{D_Matrix_linearized_Cross-Diffusion}
		\mathbf{D}\,=\,\left(\begin{array}{cc}
			\widehat{\mathcal D}_{11}&\widehat{\mathcal D}_{12} \\
			\\
			\widehat{\mathcal D}_{21} & 	\widehat{\mathcal D}_{22}
		\end{array}\right),
	\end{equation}
	with
	\begin{equation}
		\widehat{\mathcal D}_{ii}=[ c_i(\overline n_1,\,\overline n_2)]^2\,\mathcal D_i,
		\quad 
		\widehat{\mathcal D}_{ij}=- c_i(\overline n_1,\,\overline n_2)\,\mathcal D_{i}\,\overline n_i\,\lambda_{i}(\overline n_1,\,\overline n_2), \quad \mbox{for} \ \  i\ne j.
	\end{equation}
	{Similarly to the self--diffusion case, the linearized system \eqref{eq:lindiffW} with diffusion matrix  \eqref{D_Matrix_linearized_Cross-Diffusion} 
		is obtained from {the reaction--diffusion system} \eqref{SistMac2.1} 
		by using the same expansions for $c_i$ together with the ones for $\lambda_i$, i.e.
		$
		\lambda_i(n_1,\, n_2) =  \lambda_i(\overline n_1,\,\overline n_2) + {\nabla_{\bx}}  \lambda_i(\overline n_1,\,\overline n_2) \cdot \mathbf{w} + O(\mathbf{\|w\|}^2)
		$ and retaining only the linear terms. 
	}
	
	Also, the conditions leading to pattern formation, in this case, are analogous to  
	those provided in Subsection \ref{SecTuring1} and established  {in 
		\eqref{cond2.1} and \eqref{cond2.2}}.
	In this case, the conditions are
	\begin{eqnarray}
		& &{s_1(\delta, \delta_{12}, \delta_{21})  :=}  -\beta \,(\delta_{21}+1)+\delta_{12}\,(\zeta -1)^2-\delta\,\zeta +\delta+\zeta\,\nu < 0,
		\label{cond2.1c}
		\\[3mm]
		& &  \label{cond2.2c} {s_2(\delta, \delta_{12}, \delta_{21})  :=}  \left(\delta_{12}\,(\zeta -1)^2-\delta \zeta-\beta \,(\delta_{21}+1) +\delta+\zeta\,\nu\right)^2 \\[1mm]
		& & \hspace*{4cm} +\, 4 (\zeta -1) \zeta \,(\beta -\nu )\,(\delta_{12} \delta_{21}-\delta)   >0,
		\nonumber
	\end{eqnarray}
	where, as in the previous section, $\delta = \dfrac{\widehat{\mathcal D}_{22}}{\widehat{\mathcal D}_{11}}$ and  $\delta_{ij}=\dfrac{\widehat{\mathcal D}_{ij}}{\widehat{\mathcal D}_{11}}$.
	Inequalities \eqref{cond2.1c} and \eqref{cond2.2c} establish the Turing conditions of our analysis when the cross-diffusion effects
	are introduced in the dynamics, in order to observe the formation of bacterial patterns.
	{As a remark, note that the above conditions constitute a generalization of  \eqref{cond2.1}--\eqref{cond2.2}, which can indeed be written as $s_1(\delta, 0, 0) < 0$ and $s_2(\delta, 0, 0)> 0$. For any fixed {parameters} $\delta_{12}$ and $\delta_{21}$, {functions} $s_1$ and $s_2$ are polynomials in $\delta$ of, respectively, degree 1 with negative slope and degree 2 with positive leading coefficient. 
		Therefore, there will always be a point  $\overline{\delta}$ such that, if $\delta > \overline{\delta}$, {conditions  \eqref{cond2.1c} and \eqref{cond2.2c} are} satisfied. 
		In Figure \ref{fig:instability_conditions}, left panel, we plot the functions $s_1(\delta, \delta_{12}, \delta_{21})$ and $s_2(\delta, \delta_{12}, \delta_{21})$ for three different fixed choices of $ \delta_{12}$ and $\delta_{21}$. 
		In this example, we show that adding some positive  cross-diffusion in the activator equation ($\delta_{12}>0$) 
		can have stabilizing effects, as the minimal value $\overline{\delta}$,  indicated by the black small square, 
		increases with respect to the self-diffusion case of Section \ref{SecTuring1}, represented by the red small square. 
		The situation reverts with the sign of  $\delta_{12}$, as negative values have a destabilizing effect (see the small green square) or, for instance,  if one considers cross-diffusion only in the inhibitor equation, represented by the value of $\delta_{21}$,  as shown in the right panel of Figure \ref{fig:instability_conditions}.
	}
	
	{
		\begin{figure}
			\includegraphics[width = \textwidth]{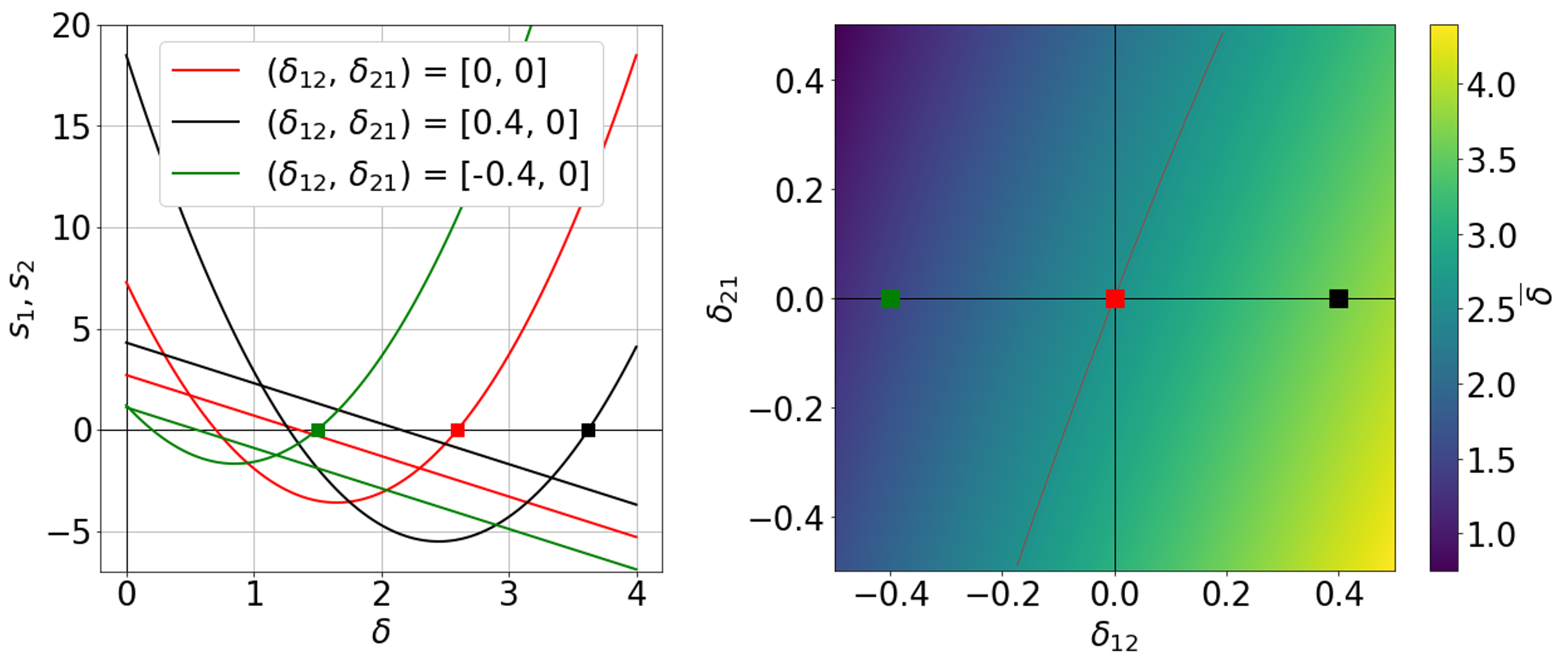}	
			\caption {{
					On the left panel: $s_1$ and $s_2$ defined in \eqref{cond2.1c}-\eqref{cond2.2c}, for three different choices of  the cross-diffusion parameters $\delta_{12}, \delta_{21}$, as functions of $\delta$.  In each case, the minimal values $\overline{\delta}$, such that  \eqref{cond2.1c}-\eqref{cond2.2c} are satisfied for all  $\delta>\overline{\delta}$, are indicated by a small square of the same colour. 
					On the right panel, the values of $\overline{\delta}$ are plotted as a function of  $\delta_{12}$ and $\delta_{21}$. In the absence of cross-diffusion (small red square),  $\overline{\delta}$ is about 2.6. Such a value is kept constant over the red line crossing the origin. 
					The values of the other parameters are 
					$\zeta = 3$,
					$\beta = 1.5$,
					$\nu = 1.4$.	
			}}
			\label{fig:instability_conditions}
		\end{figure}
	}
	In this section, we aimed to show the occurrence, in \eqref{SistMac2}, of self and cross--diffusion--driven instabilities via a linear stability analysis. However, it might be worth mentioning that, in order to better appreciate the pattern formation dynamics, one could extend the stability analysis to higher order terms, which include, for instance, the product of the gradients of $n_1$ and $n_2$, see e.g. \cite{gambino2013pattern, tulumello2014cross}.

	
	\section{Numerical simulations}
	\label{SecSim}
	
	In this section, we aim to show numerically the instability properties of the reaction-diffusion model \eqref{SistMac2.1}-\eqref{SistMac2.1_bc}. 
	Being our purposes purely illustrative, we arbitrarily fix the parameters of the model as follows:
	\begin{equation}
		\beta =1.5  ,\quad \tau = 2, \quad \nu=1.4 ,\quad {\mathcal{D}}_1=0.1,\label{MacPar}
	\end{equation}
	while we discuss the arising of Turing instability for varying parameters $\zeta$ and $\delta$. 
	Since one of the main novelties of this paper is the analytic asymptotic procedure able to consistently derive nonlinear and cross-diffusion terms from the kinetic level, in the numerical tests, we aim at investigating the effects due to turning coefficients and cellular speeds depending on both bacterial densities.
	We take, as a reference case, functions $c_i$  constantly equal to $1$ and we consider the  { turning coefficients}  
	\begin{equation}\label{Funlami1}
		\lambda_i(n_i,n_j)=0.25\left(\frac{1}{\sqrt{n_i}\,(n_i+n_j)}\right)^{\frac23}, \quad (i,j)=(1,2),(2,1),
	\end{equation}
	{so that the turning capability of bacteria decreases as the number of bacteria in a neighborhood increases.}
	In this case, conditions \eqref{cond2.1c}-\eqref{cond2.2c} rewrite as
	\begin{eqnarray}
		& &(1 - \zeta) \left( \delta + ( \zeta-1)\,\bar\lambda_1 \right) + \beta\,( \delta\,\bar\lambda_2-1) + \zeta\, \nu
		< 0,\label{cond2.1cs}
		\\[3mm]
		& &  \label{cond2.2cs} 4 \delta \, \zeta (\zeta - 1)\,( \bar\lambda_1\,\bar\lambda_2-1)\,(\beta - \nu) 
		\\[1mm]
		& & + \left((1 - \zeta)\,(\,\delta + ( \zeta-1 )\,\bar\lambda_1) + \beta\,(  \delta\,\bar\lambda_2-1) + \zeta\,\nu\right)^2
		>0,
		\nonumber
	\end{eqnarray}
	where $\bar\lambda_i=\lambda_i(\bar n_1,\bar n_2)$.
	
\begin{figure}[ht!]
		\centering
		\begin{tabular}{cc}
			\includegraphics[scale=0.3]{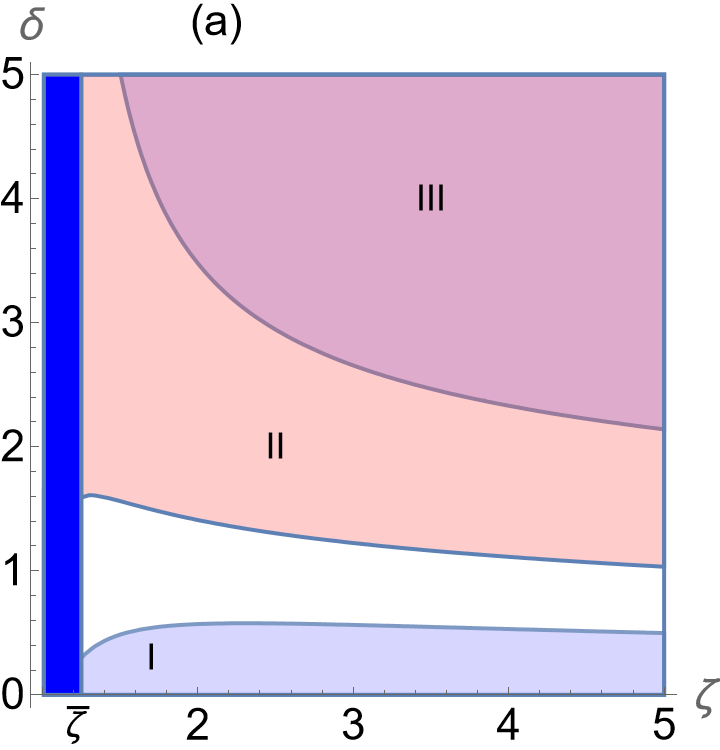} & \includegraphics[scale=0.3]{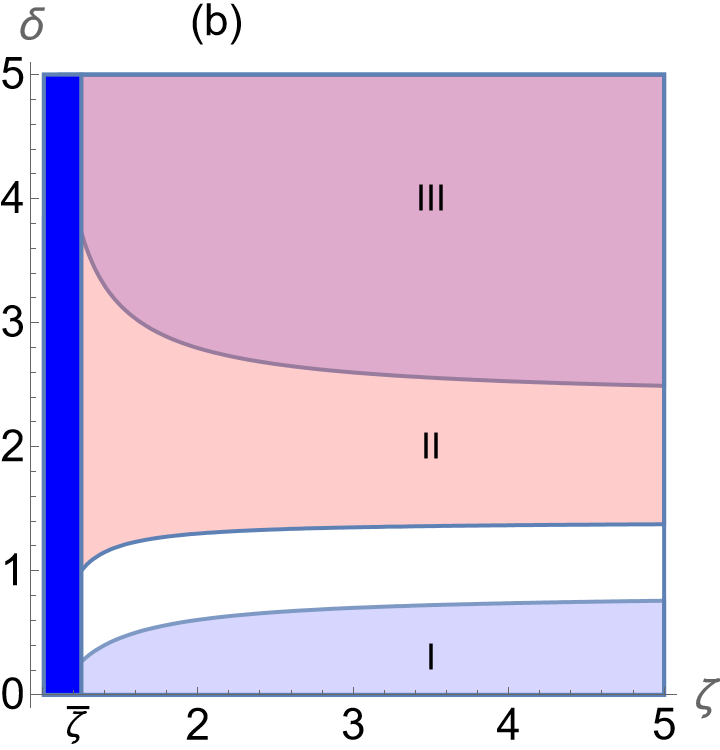}  \\
			\includegraphics[scale=0.3]{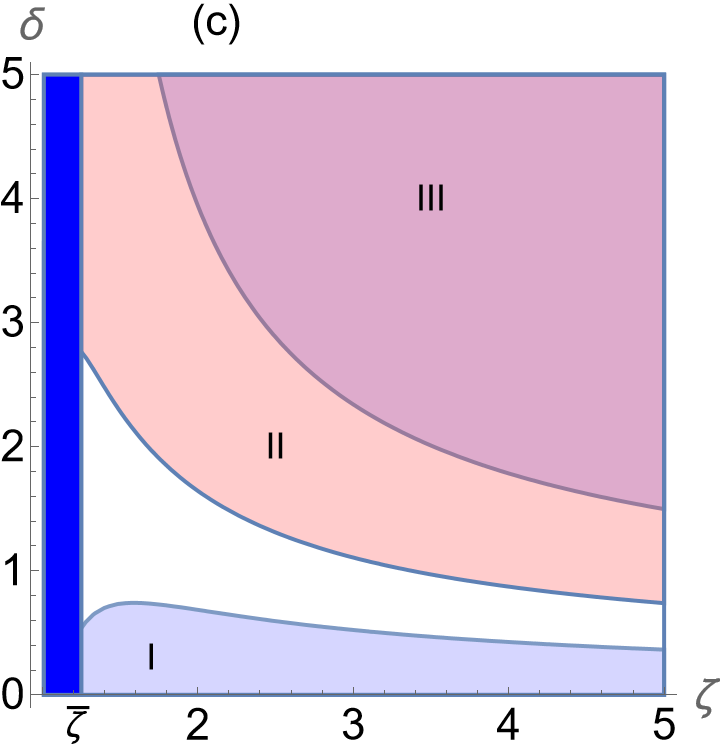} &
			\includegraphics[scale=0.3]{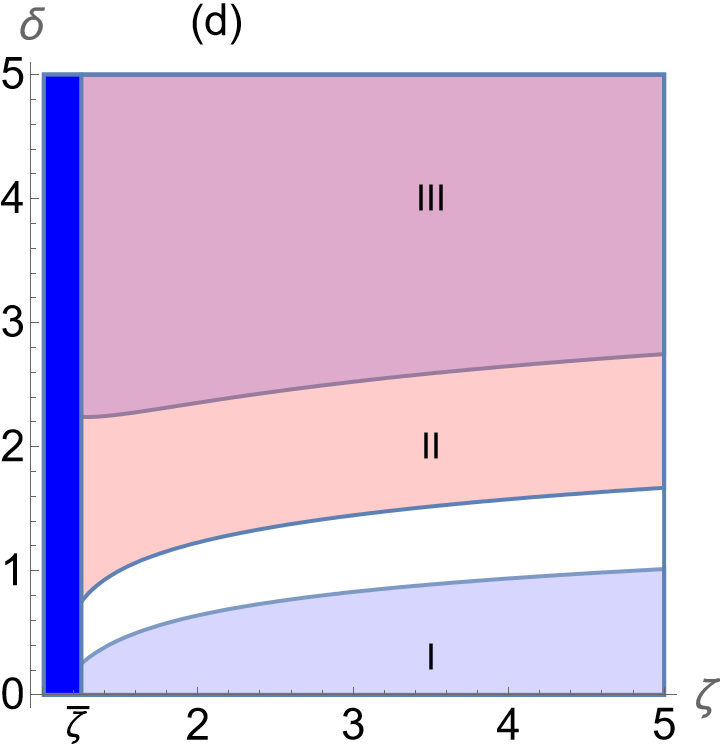} 
		\end{tabular}
		\caption {Parameters space $\zeta-\delta$ in which, holding stability condition \eqref{zetabar}, conditions for Turing instability and formation of patterns are satisfied (region III in each panel), taking functions $c_i$ constantly equal to $1$  (panels (a), (b), (d)) or as in \eqref{Funci1} (panel (c)) and functions $\lambda_i$ as in \eqref{Funlami1}  (panels (a), (c)), null (panel (b)) or as in \eqref{Funlami2}  (panel (d)). Other parameters are taken as in \eqref{MacPar}. 	
		}
		\label{BifDiags}       
	\end{figure}

	We observe that the present choice for functions $\lambda_i$ corresponds to an attractive behavior of the two strains. In other words, bacteria of one population tend to reach the individuals of the other one, and this enhances the prevalence of the cooperative process. 
	We report the bifurcation diagram for this case in Figure \ref{BifDiags}, panel (a), where we show the values for which, holding conditions \eqref{zetabar}, conditions \eqref{cond2.1cs} (region I), \eqref{cond2.2cs} (region II), or both  (region III) are satisfied.
	
\begin{figure}
		\includegraphics[width = \textwidth]{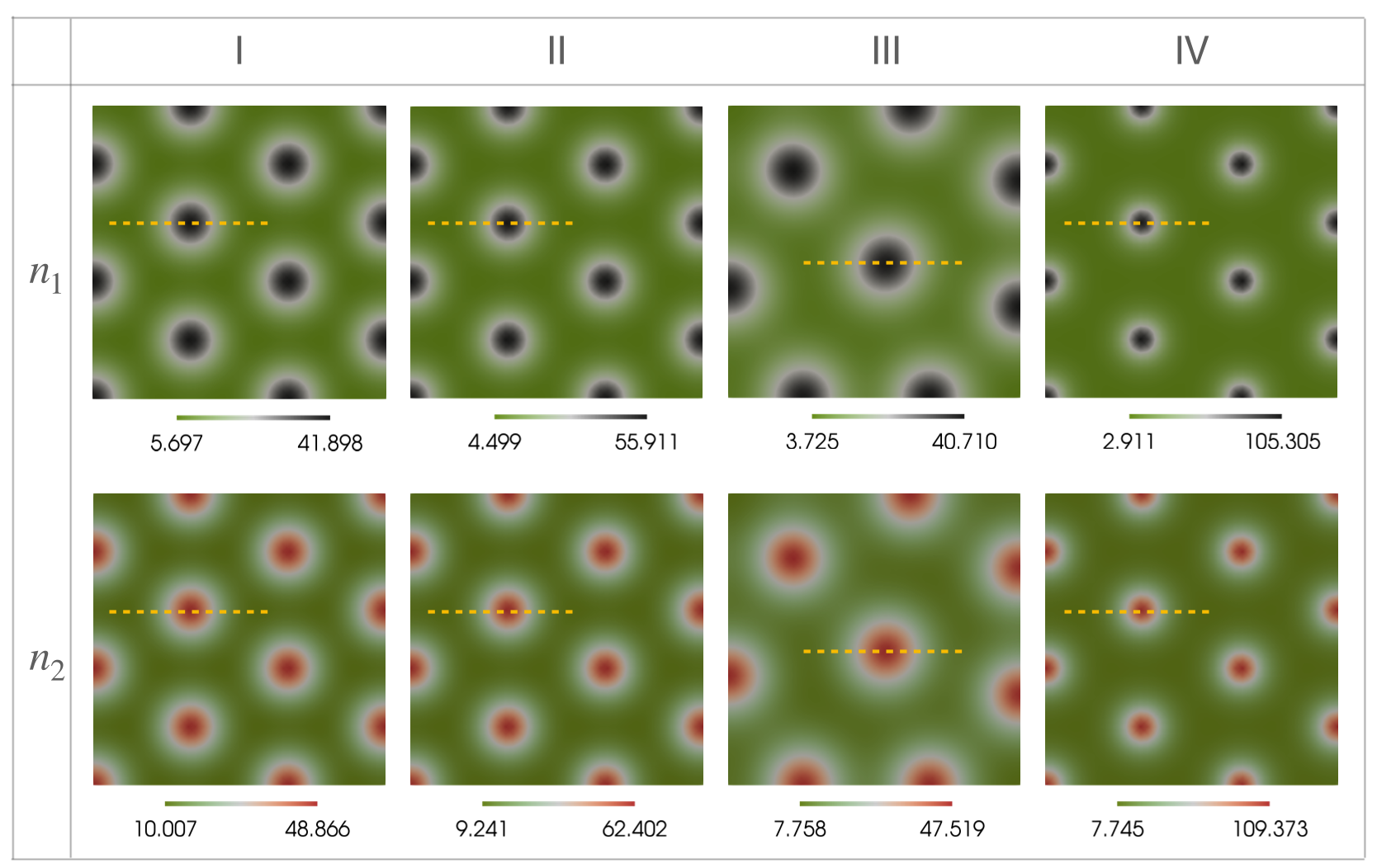}
		\caption{
			A collection of four different numerical solutions to \eqref{SistMac2.1}-\eqref{SistMac2.1_bc}, for different choices of the diagonal and cross-diffusion coefficients,  at the final time $T=1000$. 
			In the first column,  I,  the solutions $n_1$ and $n_2$ for $c_1=c_2=1$ and cross-diffusion terms $\lambda_1,\lambda_2$ as in  \eqref{Funlami1} are shown.  In the second column,  II,  $c_1=c_2=1$ and $\lambda_1 = \lambda_2 = 0$ (i.e.  the system reduces to the self-diffusion case \eqref{SistMac2}).  In III, the diagonal diffusion has terms $c_1, c_2$ as in \eqref{Funci1} and cross-diffusion as in \eqref{Funlami1}.  In the last case, IV, $c_1=c_2=1$ and cross-diffusion functions as in \eqref{Funlami2}.  The yellow dashed segments over the spots refer to Figure \ref{fig:numerical simulations_profiles}, where the profiles of the spots along such segments are compared for the 4 cases.
			The remaining parameters are
			$\zeta = 3$,
			$\beta = 1.5$,
			$\tau = 2$,
			$\nu = 1.4$,
			$D_1 = 0.1$, 
			$\delta_{22} = 2.7$.
		}
		\label{fig:numerical simulations}
	\end{figure}

	We take $\zeta=3$ and $\delta=2.7$ and perform numerical simulations for this case on a square domain ${\Gamma_{\bx}} = [0,\pi]\times[0,\pi]$.  The numerical method, which combines finite elements in space and finite differences in time,  is described in Appendix \ref{Appendix:numerical_method}.  Results showing aggregation of bacteria forming spots on the surface are shown in Figure \ref{fig:numerical simulations}, column I.
	\begin{figure}
		\includegraphics[width = \textwidth]{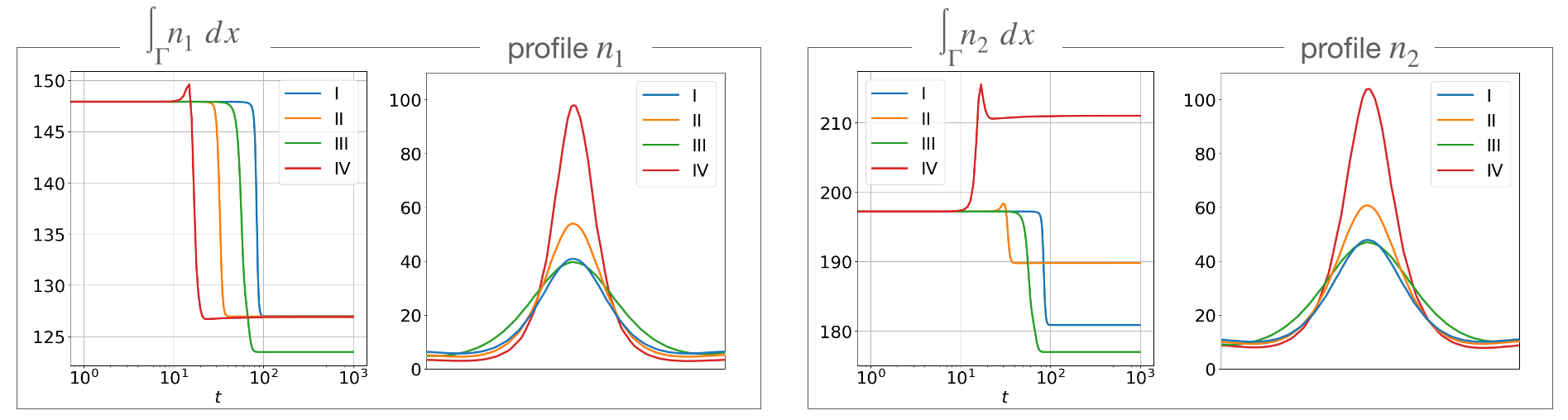}
		\caption{
			Some details of the solutions $n_1$ and $n_2$ are reported in Figure \ref{fig:numerical simulations}.  The first box refers to the function $n_1$ for the cases I--IV.  In the first image, we plot its integral over time, in the second image, the profile of the spots,  along the cuts highlighted in Figure \ref{fig:numerical simulations} with yellow dashed segments. 
			In the same way,  the second box relates to the function $n_2$.
		}
		\label{fig:numerical simulations_profiles}
	\end{figure}
	
	Now we compare this case to the one in which no cross-diffusion is included, namely $\lambda_1=\lambda_2=0$. In this framework, the stability conditions are \eqref{cond2.1}-\eqref{cond2.2}, and the bifurcation diagram relevant to this case is reported in Figure \ref{BifDiags}, panel (b). Taking again $\zeta=3$ and $\delta=2.7$, we show the behavior of system \eqref{SistMac2} in Figure \ref{fig:numerical simulations} (II).
	
	Comparing these results with the previous ones, we observe that the configuration in space of the two strains is analogous.  However, the presence of the cross-diffusion terms in case II seems to have an inhibition effect on the overall growth of $n_2$,  as reported by the corresponding integral in Figure \ref{fig:numerical simulations_profiles}. This is due to the fact that, being the density of both the strains higher at the center of each spot, the cooperation dynamics result in a growth for $n_2$. On the other hand, because of inter-specific competition, the total mass of $C_1$ reaches the same value at the steady state of both cases I and II.

	In the third case, instead, we keep the cross-diffusion term as in \eqref{Funlami1} and we take functions $c_i$ as
	\begin{equation}\label{Funci1}
		c_i(n_i,n_j)=1+0.5\,\left(\frac{n_i}{n_i+n_j}\right)^{\frac23}, \quad
		(i,j)=(1,2),(2,1).
	\end{equation} 
	This choice leads to a diffusion term incremented by a term that depends on the fraction of the strain density over the total one, i.e., the higher the fraction of concentration of their strain is, the more individuals will tend to diffuse. 
	This implies that conditions \eqref{cond2.1c}-\eqref{cond2.2c} now become
	\begin{eqnarray}
		& &(1 - \zeta) \left( \bar c\, \delta + ( \zeta-1)\,\bar\lambda_1 \right) + \beta\,( \delta\,\bar\lambda_2-1) + \zeta\, \nu
		< 0,\label{cond2.1css}
		\\[3mm]
		& &  \label{cond2.2css} 4 \delta \, \zeta (\zeta - 1)\,( \bar\lambda_1\,\bar\lambda_2-\bar c)\,(\beta - \nu) 
		\\[1mm]
		& & + \left((1 - \zeta)\,(\bar c\,\delta + ( \zeta-1 )\,\bar\lambda_1) + \beta\,(  \delta\,\bar\lambda_2-1) + \zeta\,\nu\right)^2
		>0,
		\nonumber
	\end{eqnarray}
	where, in this case, $\bar\lambda_1=\dfrac{\lambda_1(\bar n_1,\bar n_2)}{c_1(\bar n_1,\bar n_2)}$, $\bar\lambda_2=\dfrac{\lambda_2(\bar n_1,\bar n_2)\,c_2(\bar n_1,\bar n_2)}{(c_1(\bar n_1,\bar n_2))^2}$ and $\bar c=\left(\dfrac{c_2(\bar n_1,\bar n_2)}{c_1(\bar n_1,\bar n_2)}\right)^2$. 
	The bifurcation diagram relevant to this case is reported in Figure \ref{BifDiags}, panel (c). We report numerical simulations for this third case in Figure \ref{fig:numerical simulations} (case III), taking $\zeta=3$ as before and $\delta=2.41$. Comparing the results with the reference case of constant diffusion, we may observe that the bigger diffusion coefficient leads to a major spreading of cells, resulting in a lower number of spots. Moreover, both the total masses of bacteria strains $C_1$ and $C_2$ are reduced with respect to the other cases (see Figure \ref{fig:numerical simulations_profiles}), due to a lack of cooperative behavior.
	
	As last case, we take again functions $c_i$ constantly equal to $1$ and we consider cross-diffusion terms as follows
	\begin{equation}\label{Funlami2}
		\lambda_i(n_i,n_j)=-0.25\left(\frac{1}{\sqrt{n_i}\,(n_i+n_j)}\right)^{\frac23}, \quad (i,j)=(1,2),(2,1).
	\end{equation}
	Specifically, these cross-diffusion terms have opposite signs with respect to the reference case, meaning that the two bacterial populations undergo repulsive dynamics, avoiding one another. 
	The bifurcation diagram pertinent to this scenario is shown in Figure \ref{BifDiags}, panel (d). We present numerical simulations for this further case in Figure \ref{fig:numerical simulations} (case IV), with $\zeta=3$ and $\delta=2.7$ again. We may observe that, in this setting, the number of spots visible on the domain is the same as in the attractive case. The outstanding difference is the fact that the concentration of bacteria is considerably higher at the center of each spot, while they tend to be less present in areas between spots, where the repulsion appears to be prevalent. Lastly, among the four cases, this scenario is the most proliferative one for the species $n_2$, which increases with respect to its initial value, while the mass of $C_1$ is substantially the same as in cases I and II. 
	
	
	\section{Concluding remarks and perspectives}
	\label{SecConc}
	
	In this work, we have proposed a mathematical model capable of comprehensively describing a biological system, 
	in which individuals from different species may interact among themselves and with their environment. 
	Our approach aims to capture the dynamics at both a mesoscopic level, 
	as a result of straightforward interactions among the individuals,
	and a macroscopic level, where the observable phenomena are directly linked to the mesoscopic description. 
	To do so, we employed the kinetic theory of active particles, which enabled us to describe
	all possible interactions through integral operators.
	More specifically, these interactions can be of a conservative type, 
	where individuals may only change their internal state or velocity, 
	and of nonconservative type, where the interplay leads to proliferative or destructive phenomena. 
	Furthermore, we accounted for cooperative interactions among species, 
	as well as for intraspecific and interspecific competitive interactions.
	We have also assumed that the speed of the individuals depends on their internal state. 
	In order to obtain a macroscopic system capable of describing collective trends, 
	such as the formation of spatial patterns, the mesoscopic description has been performed, taking into account the specific timescale for each interactive process. 
	In our case, we have assumed that the conservative interactions of the individuals with the host constitute the dominant process, occurring at a faster time scale. 
	Moreover, we have considered that changes in the direction of velocity,
	{influenced} by the macroscopic densities of the other species,
	occurs at a time scale faster than the one for nonconservative phenomena. 
	After stating an analytical result for the interactions with the host, we have derived macroscopic reaction-cross diffusion equations for the global densities of the species,
	through a hydrodynamic limit. 
	
	The procedure outlined has been applied to the biological setting of two bacterial strains on a leaf surface. 
	The mesoscopic description has allowed us to focus on the cooperative and competitive interplay between the two strains, 
	as well as the nourishing function of the leaf cells. 
	We have then derived a macroscopic system of two partial differential equations,
	where the self-diffusion and the cross-diffusion coefficients can be assumed as functions of the macroscopic densities themselves. 
	
	Performing numerical simulations, we have observed how the choices of these functions may influence the effects 
	of the cooperation or the competition, leading to either enhanced or diminished co-aggregation. 
	{
		Despite its consistency and applicability to a real biological system, our model presents some limitations that could be overcome in future research.
		Indeed, the macroscopic equations we derive could be generalized to more complex multiscale frameworks.
		Moreover, our asymptotic procedure relies on specific scaling assumptions, which could potentially be relaxed in future work to capture additional relevant biological time scales.
		Similarly, our approach could be further unified along the lines of the upscaling methods presented in other works.
		The present model could also be extended to include more detailed cell-level mechanisms, in order to bring it closer to the sophistication of current multiscale approaches while retaining its focus on ecological interactions.}
	
	{In addition, the development of models that explicitly account for the spatial heterogeneity of the host environment and its interactions with bacterial populations represents a promising avenue for future research. For instance, incorporating the anisotropic structure of the leaf surface and allowing the distribution functions of both host and nourishing cells 
		to vary spatially in response to bacterial activity could enable the model to capture key phenomena such as localized changes in the host medium.}
	Moreover, a further investigation of the obtained macroscopic systems, including rigorous results on the global existence and uniqueness of 
	a nonnegative solution 
	(as those provided in \cite{eckardt2024mathematical}),
	may lead to a more refined study of the problem.
	
	{Finally, future perspectives include performing a weakly nonlinear analysis of the macroscopic system in order to fully characterize the types of patterns that can arise and their stability properties.}
	
	\section*{Acknowledgments}
	
	The work of RT and MB was performed in the frame of activities sponsored by the Italian National Group of Mathematical Physics (GNFM-INdAM) and by the University of Parma (Italy). 

\section*{Funding Declaration}

 The work of RT and MB was carried out within the activities sponsored by the Italian National Group of Mathematical Physics (GNFM–INdAM) and by the University of Parma (Italy). RT is a postdoctoral fellow supported by the National Institute of Advanced Mathematics (INdAM), Italy.

MB acknowledges support from the project PRIN 2022 PNRR “Mathematical Modelling for a Sustainable Circular Economy in Ecosystems” (project code P2022PSMT7, CUP D53D23018960001), funded by the European Union – NextGenerationEU PNRR-M4C2-I 1.1 and by MUR–Italian Ministry of Universities and Research.

MB and RT also acknowledge support from the University of Parma through the Bando di Ateneo 2022 per la ricerca, cofunded by MUR–Italian Ministry of Universities and Research (D.M. 737/2021 – PNR – PNRR – NextGenerationEU), project “Collective and Self-Organised Dynamics: Kinetic and Network Approaches”.

The work of DC, AJS, and RT was supported by COST Action CA18232, by Project UID/00013/25: Centro de Matemática da Universidade do Minho (CMAT/UM), and by Portuguese national funds (OE) through the project FCT/MCTES PTDC/03091/2022, “Mathematical Modelling of Multi-scale Control Systems: applications to human diseases – CoSysM3” (\url{https://doi.org/10.54499/2022.03091.PTDC}
).

\section*{CRediT Author Statement}

All authors have contributed equally to this work. The following roles were performed jointly by all authors: Conceptualization; Data curation; Formal analysis; Funding acquisition; Investigation; Methodology; Project administration; Resources; Software; Supervision; Validation; Visualization; Writing – original draft; Writing – review \& editing.
Each author has read and approved the final manuscript.

	
	\appendix
	
	\section{Appendix -- The numerical method}
	\label{Appendix:numerical_method}
	
	In this section, we briefly present the method we use to solve system \eqref{SistMac2.1}, which, for convenience, we report here  in the following form:
	\begin{equation}
		\label{numerical:SistMac2.1}
		\begin{aligned}
			\frac{\pa n_1}{\pa t} \!&=\! \, \,c_1\,{\mathcal{D}}_1\,\nabla_{\bx}\cdot \left(c_1\,\,\nabla_{\bx}\, n_1\,
			-\lambda_{1}\,n_1\,\nabla_{\bx} n_2\right)+ f_1(n_1,n_2) n_1 \, , 
			\\[2mm]
			\frac{\pa n_2}{\pa t} \!&=\!  \,c_2\,{\mathcal{D}}_2\,\nabla_{\bx}\cdot \left(c_2\,\nabla_{\bx}\, n_2\,
			-\lambda_{2}\,n_2\,\nabla_{\bx} n_1\right)+f_2(n_1,n_2) n_2 \, ,
		\end{aligned}
	\end{equation}
	with 
	$$
	f_1(n_1,n_2) = \frac{\zeta\,n_1\,n_2}{n_1+\beta\,n_2} - \,n_2, 
	\quad \text{and} \quad 
	f_2(n_1,n_2) = \frac{\zeta\,\beta\, n_1\,n_2}{n_1+\beta\,n_2} -\tau-\nu\,n_2.
	$$ 
	The numerical method is based on the application of the finite element method in space and finite differences in time. Let $X = <\varphi_1(x), \ldots, \varphi_{N_x}(x)>$ be the finite element space defined on a triangulation ${\color{magenta}{{\Gamma_{\bx}}}}^h$ of the domain ${\Gamma_{\bx}}$, where we look for numerical solutions $\widetilde{n_1} = \sum_{j =1}^{N_x} a_{j}(t)\varphi_j(x)$ and $\widetilde{n_2} = \sum_{j=1}^{N_x} b_{j}(t)\varphi_j(x)$ to approximate, respectively, the solutions $n_1$ and $n_2$ of \eqref{numerical:SistMac2.1}.
	By defining the vector of coefficients $\textbf{n}(t) := \left(a_1, \ldots, a_{N_x}, b_1, \ldots, b_{N_x} \right)^t$, \eqref{numerical:SistMac2.1} is discretized in space by the ODE system
	\begin{equation}
		\label{numerical:ODE}
		\textbf{n}' 
		+ K(\textbf{n}) \textbf{n} 
		= F(\textbf{n}) \textbf{n},    
	\end{equation}
	where 
	$$
	K = \begin{pmatrix}
		K_{1,1}& 
		K_{1,2} \\[1mm]
		K_{2,1} & 
		K_{2,2}
	\end{pmatrix}
	$$
	is a block matrix, with 
	\begin{align*}
		& K_{1,1} = \left(
		\int_{{\color{magenta}{{\Gamma_{\bx}}}}^h}\mathcal{D}_{1} c_1\nabla \varphi_j \cdot \nabla (c_1\varphi_i)\;dx
		\right)_{1\leq i,j\leq N_x}
		\\
		&
		K_{1,2} = \left(
		- \int_{{\color{magenta}{{\Gamma_{\bx}}}}^h}\mathcal{D}_{1} \lambda_1 
		\nabla \varphi_j \cdot \nabla (c_1\varphi_i)\;dx
		\right)_{1\leq i,j\leq N_x}
		\\
		& 
		K_{2,1} = \left(
		- \int_{{\color{magenta}{{\Gamma_{\bx}}}}^h}\mathcal{D}_{2} \lambda_2 
		\nabla \varphi_j \cdot \nabla (c_2\varphi_i)\;dx
		\right)_{1\leq i,j\leq N_x}
		\\
		& 
		K_{2,2} = \left(
		\int_{{\color{magenta}{{\Gamma_{\bx}}}}^h}\mathcal{D}_{2} c_2\nabla \varphi_j \cdot \nabla (c_2\varphi_i)\;dx
		\right)_{1\leq i,j\leq N_x}
	\end{align*}
	and 
	$$
	F = \begin{pmatrix}
		F_{1}& 
		0 \\
		0 & 
		F_{2}
	\end{pmatrix}
	$$ is a block matrix with 
	$$
	F_1 = \left(
	\int_{{\color{magenta}{{\Gamma_{\bx}}}}^h} 
	f_1 \varphi_j \varphi_i
	\;dx
	\right)_{1\leq i,j\leq N_x}
	\text{and}  \quad
	F_2 = \left(
	\int_{{\color{magenta}{{\Gamma_{\bx}}}}^h} 
	f_2 \varphi_j \varphi_i
	\;dx
	\right)_{1\leq i,j\leq N_x}.
	$$
	We remark that, since the functions $c_1, c_2, \lambda_1, \lambda_2, f_1, f_2$ depend on $\widetilde{n_1}$ and $\widetilde{n_2}$, both $K$ and $F$ depend on $\textbf{n}$.
	Next, introduce a timestep $\Delta t>0$ and, by applying a backward Euler scheme, we get the following fully implicit approximation of the ODE system \eqref{numerical:ODE}
	\begin{equation}
		\label{numerical:discrete}
		\frac{\textbf{n}^{k+1}}{\Delta t} + K(\textbf{n}^{k+1}) \textbf{n}^{k+1} - F(\textbf{n}^{k+1}) \textbf{n}^{k+1} = \frac{\textbf{n}^{k}}{\Delta t}, 
	\end{equation}
	which is a nonlinear system. 
	In order to solve the nonlinearity, we apply a fixed point iteration method with tolerance $\varepsilon>0$ (see e.g  \cite{madzvamuse2014fully}):
	at every time-step $t_{k+1} = (k+1)\Delta t$, we set $\textbf{p}^0 = \textbf{n}^{k}$ and, 
	for $j\geq 0$, we iteratively solve the linear system
	\begin{equation}\label{numerical:picard}
		\frac{\textbf{p}^{j+1}}{\Delta t} + K(\textbf{p}^{j}) \textbf{p}^{j+1} - F(\textbf{p}^{j}) \textbf{p}^{j+1} = \frac{\textbf{n}^{k}}{\Delta t}.    
	\end{equation}
	If there exists $\overline{j}$ such that $\| \textbf{p}^{\overline{j}+1} - \textbf{p}^{\overline{j}} \|_2 < \varepsilon$, we stop the iterations and set $\textbf{n}^{k+1} := \textbf{p}^{\overline{j}+1}$.
	
	The numerical scheme is implemented in Python 3, and it is solved by using FEniCS \cite{logg2012automated} (dolfin version 2019.1.0).  The mesh of the spatial domain ${\color{magenta}{{\Gamma_{\bx}}}} = [0,\pi]\times[0,\pi]$ is created by a $40\times40$ homogeneous partition, where each small square is divided into 4 equilateral and identical triangles, for a total of 6400 triangles.  We consider the function space $X$ constituted by continuous and piecewise linear functions.  For the temporal discretization, we set $\Delta t = 0.01$ and $\varepsilon=0.01$.
	
	In all of the simulations, the initial condition $(\overline{n}_1, \overline{n}_2)$, constituted by the homogeneous steady state \eqref{Equi},  
	is perturbed at each vertex of the mesh by the addition of a random value sampled from a Gaussian distribution 
	with mean 0 and variance of, respectively, $\overline{n}_1/100$ and $\overline{n}_2/100$. 
	For all of the cases presented in this work, the perturbation has the same seed(1) from the Python library \emph{random}, 
	which makes the initial profiles equivalent for all the simulations. 
	

	\bibliographystyle{abbrv}
	\bibliography{biblio_new}

\begin{thebibliography}{10}

\bibitem{abbas2010existence}
S.~Abbas, M.~Banerjee, and N.~Hungerb{\"u}hler.
\newblock Existence, uniqueness and stability analysis of allelopathic
  stimulatory phytoplankton model.
\newblock {\em J. Math. Anal. Appl.}, 367(1):249--259, 2010.

\bibitem{alt1980biased}
W.~Alt.
\newblock Biased random walk models for chemotaxis and related diffusion
  approximations.
\newblock {\em J. Math. Biology}, 9:147--177, 1980.

\bibitem{Anwasia-etal-2017}
B.~Anwasia, P.~Gonçalves, and A.~J. Soares.
\newblock From the simple reacting sphere kinetic model to the
  reaction-diffusion system of {M}axwell-{S}tefan type.
\newblock {\em Commun. Math. Sci.}, 17(2):507--538, 2019.

\bibitem{bardos2016simultaneous}
C.~Bardos and H.~Hutridurga.
\newblock Simultaneous diffusion and homogenization asymptotic for the linear
  {B}oltzmann equation.
\newblock {\em Asympt. Anal.}, 100(1-2):111--130, 2016.

\bibitem{Bellomo-Belloquid}
N.~Bellomo and A.~Belloquid.
\newblock From a class of kinetic models to the macroscopic equations for
  multicellular systems in biology.
\newblock {\em Discrete Contin. Dyn. Syst., Ser. B.}, 4(1):59--80, 2004.

\bibitem{survey}
N.~Bellomo, D.~Burini, G.~Dosi, L.~Gibelli, D.~Knopoff, N.~Outada, P.~Terna,
  and M.~E. Virgillito.
\newblock What is life? {A} perspective of the mathematical kinetic theory of
  active particles.
\newblock {\em Math. Models Methods Appl. Sci.}, 31(09):1821--1866, 2021.

\bibitem{bellomo1994dynamics}
N.~Bellomo and G.~Forni.
\newblock Dynamics of tumor interaction with the host immune system.
\newblock {\em Math. Comput. Modelling}, 20(1):107--122, 1994.

\bibitem{bendahmane2024mathematical}
M.~Bendahmane, F.~Karami, D.~Meskine, J.~Tagoudjeu, and M.~Zagour.
\newblock Mathematical analysis and multiscale derivation of a nonlinear
  predator--prey cross-diffusion--fluid system with two chemicals.
\newblock {\em Commun. Nonlinear Sci. Numer. Simul.}, 136:108090, 2024.

\bibitem{bertotti2023modelling}
M.~L. Bertotti, B.~Carbonaro, and M.~Menale.
\newblock Modelling a market society with stochastically varying money exchange
  frequencies.
\newblock {\em Symmetry}, 15(9):1751, 2023.

\bibitem{bisi2006reactive}
M.~Bisi and L.~Desvillettes.
\newblock From reactive {B}oltzmann equations to reaction-diffusion systems.
\newblock {\em J. Stat. Phys.}, 124:881--912, 2006.

\bibitem{bisi2025derivation}
M.~Bisi, M.~Groppi, G.~Martal{\`o}, and R.~Travaglini.
\newblock Derivation from kinetic theory and 2-{D} pattern analysis of
  chemotaxis models for {M}ultiple {S}clerosis.
\newblock {\em J. Math. Biol.}, 91(43), 2025.

\bibitem{bisi2022reaction}
M.~Bisi and R.~Travaglini.
\newblock Reaction-diffusion equations derived from kinetic models and their
  {T}uring instability.
\newblock {\em Commun. Math. Sci.}, 20(3):763--801, 2022.

\bibitem{brewer1991functional}
C.~Brewer, W.~Smith, and T.~Vogelmann.
\newblock Functional interaction between leaf trichomes, leaf wettability and
  the optical properties of water droplets.
\newblock {\em Plant Cell Environ.}, 14(9):955--962, 1991.

\bibitem{brezis2010functional}
H.~Brezis.
\newblock {\em Functional Analysis, Sobolev Spaces and Partial Differential
  Equations}.
\newblock Universitext. Springer, 2010.

\bibitem{burch2014hygroscopic}
A.~Y. Burch, V.~Zeisler, K.~Yokota, L.~Schreiber, and S.~E. Lindow.
\newblock The hygroscopic biosurfactant syringafactin produced by {P}seudomonas
  syringae enhances fitness on leaf surfaces during fluctuating humidity.
\newblock {\em Environ. Microbiol.}, 16(7):2086--2098, 2014.

\bibitem{burini2019multiscale}
D.~Burini and N.~Chouhad.
\newblock A multiscale view of nonlinear diffusion in biology: {F}rom cells to
  tissues.
\newblock {\em Math. Models Methods Appl. Sci.}, 29(04):791--823, 2019.

\bibitem{burkhardt1999measurements}
J.~Burkhardt, H.~Kaiser, H.~Goldbach, and L.~Kappen.
\newblock Measurements of electrical leaf surface conductance reveal
  re-condensation of transpired water vapour on leaf surfaces.
\newblock {\em Plant Cell Environ.}, 22(2):189--196, 1999.

\bibitem{canizo2018rate}
J.~A. Ca{\~n}izo, A.~Einav, and B.~Lods.
\newblock On the rate of convergence to equilibrium for the linear {B}oltzmann
  equation with soft potentials.
\newblock {\em J. Math. Anal. Appl.}, 462(1):801--839, 2018.

\bibitem{conte2023mathematical}
M.~Conte, Y.~Dzierma, S.~Knobe, and C.~Surulescu.
\newblock Mathematical modeling of glioma invasion and therapy approaches via
  kinetic theory of active particles.
\newblock {\em Math. Models Methods Appl. Sci.}, 33(05):1009--1051, 2023.

\bibitem{cordier2005kinetic}
S.~Cordier, L.~Pareschi, and G.~Toscani.
\newblock On a kinetic model for a simple market economy.
\newblock {\em J. Stat. Phys.}, 120:253--277, 2005.

\bibitem{della2021sir}
R.~Della~Marca, N.~Loy, and A.~Tosin.
\newblock An {S}{I}{R}-like kinetic model tracking individuals' viral load.
\newblock {\em Netw. Heterog. Media}, 17(3):467--494, 2022.

\bibitem{della2022mathematical}
R.~Della~Marca, M.~P. Machado~Ramos, C.~Ribeiro, and A.~J. Soares.
\newblock Mathematical modelling of oscillating patterns for chronic autoimmune
  diseases.
\newblock {\em Math. Methods Appl. Sci.}, 45(11):7144--7161, 2022.

\bibitem{dietrich2022multiscale}
A.~Dietrich, N.~Kolbe, N.~Sfakianakis, and C.~Surulescu.
\newblock Multiscale modeling of glioma invasion: from receptor binding to
  flux-limited macroscopic {P}{D}{E}s.
\newblock {\em Multiscale Model. Simul.}, 20(2):685--713, 2022.

\bibitem{eckardt2024mathematical}
M.~Eckardt and C.~Surulescu.
\newblock On a mathematical model for cancer invasion with repellent p{H}-taxis
  and nonlocal intraspecific interaction.
\newblock {\em Z. Angew. Math. Phys.}, 75(2):41, 2024.

\bibitem{esser2015spatial}
D.~S. Esser, J.~H. Leveau, K.~M. Meyer, and K.~Wiegand.
\newblock Spatial scales of interactions among bacteria and between bacteria
  and the leaf surface.
\newblock {\em FEMS Microbiol. Ecol.}, 91(3), 2014.

\bibitem{franklin2007statistical}
R.~B. Franklin and A.~L. Mills.
\newblock {\em The Spatial Distribution of Microbes in the Environment}.
\newblock Springer, 2007.

\bibitem{gambino2013pattern}
G.~Gambino, M.~C. Lombardo, and M.~Sammartino.
\newblock Pattern formation driven by cross-diffusion in a 2{D} domain.
\newblock {\em Nonlinear Anal. Real World Appl.}, 14(3):1755--1779, 2013.

\bibitem{harshey2003bacterial}
R.~Harshey.
\newblock Bacterial motility on a surface: many ways to a common goal.
\newblock {\em Ann. Rev. in Microbiol.}, 57(1):249--273, 2003.

\bibitem{holscher2017sliding}
T.~H\"olscher and {\'A}.~T. Kov{\'a}cs.
\newblock Sliding on the surface: bacterial spreading without an active motor.
\newblock {\em Environ. Microbiol.}, 19(7):2537--2545, 2017.

\bibitem{krimm2005epiphytic}
U.~Krimm, D.~Abanda-Nkpwatt, W.~Schwab, and L.~Schreiber.
\newblock Epiphytic microorganisms on strawberry plants ({F}ragaria ananassa
  cv. {E}lsanta): identification of bacterial isolates and analysis of their
  interaction with leaf surfaces.
\newblock {\em FEMS Microbiol. Ecol.}, 53(3):483--492, 2005.

\bibitem{kunzler2024hitching}
M.~Kunzler, R.~O. Schlechter, L.~Schreiber, and M.~N.~P. Remus-Emsermann.
\newblock Hitching a ride in the phyllosphere: {S}urfactant production of
  {P}seudomonas spp. causes co-swarming of {P}antoea eucalypti 299{R}.
\newblock {\em Microb. Ecol.}, 87(1):62, 2024.

\bibitem{lachowicz2002microscopic}
M.~Lachowicz.
\newblock From microscopic to macroscopic description for generalized kinetic
  models.
\newblock {\em Math. Models Methods Appl. Sci.}, 12(07):985--1005, 2002.

\bibitem{logg2012automated}
A.~Logg, K.-A. Mardal, and G.~Wells.
\newblock {\em Automated solution of differential equations by the finite
  element method: The FEniCS book}, volume~84.
\newblock Springer Science \& Business Media, 2012.

\bibitem{ma1996mathematical}
Z.~Ma.
\newblock {\em Mathematical Modeling and Research on the Population Ecology}.
\newblock Anhui Education Publishing House, Hefei, 1996.

\bibitem{ramos2019kinetic}
M.~P. Machado~Ramos, C.~Ribeiro, and A.~J. Soares.
\newblock A kinetic model of {T} cell autoreactivity in autoimmune diseases.
\newblock {\em J. Math. Biology}, 79(6-7):2005--2031, 2019.

\bibitem{madzvamuse2014fully}
A.~Madzvamuse and A.~H. Chung.
\newblock Fully implicit time-stepping schemes and non-linear solvers for
  systems of reaction-diffusion equations.
\newblock {\em Appl. Math. Comput.}, 244:361--374, 2014.

\bibitem{BMT}
G.~Martalò, S.~Boccelli, and R.~Travaglini.
\newblock Turing instability and 2-{D} pattern formation in reaction-diffusion
  systems derived from kinetic theory.
\newblock {\em arXiv preprint arXiv:2509.20268}, 2025.

\bibitem{monier2003differential}
J.-M. Monier and S.~Lindow.
\newblock Differential survival of solitary and aggregated bacterial cells
  promotes aggregate formation on leaf surfaces.
\newblock {\em Proceedings of the National Academy of Sciences},
  100(26):15977--15982, 2003.

\bibitem{monier2004frequency}
J.-M. Monier and S.~Lindow.
\newblock Frequency, size, and localization of bacterial aggregates on bean
  leaf surfaces.
\newblock {\em Appl. Environ. Microbiol.}, 70(1):346--355, 2004.

\bibitem{morris1997methods}
C.~E. Morris, J.~Monier, and M.~Jacques.
\newblock Methods for observing microbial biofilms directly on leaf surfaces
  and recovering them for isolation of culturable microorganisms.
\newblock {\em Appl. Environ. Microbiol.}, 63(4):1570--1576, 1997.

\bibitem{mu2023hopf}
Y.~Mu and W.-C. Lo.
\newblock Hopf and {T}uring bifurcation for a competition and cooperation
  system with spatial diffusion effect.
\newblock {\em J. Comput. Appl. Math.}, 422:114924, 2023.

\bibitem{murray2003mathematical}
J.~Murray.
\newblock {\em Mathematical Biology II: Spatial Models and Biomedical
  Applications}, volume~18 of {\em Interdisciplinary Applied Mathematics}.
\newblock Springer, New York, 3rd edition, 2003.

\bibitem{oliveira2024reaction}
J.~M. Oliveira and R.~Travaglini.
\newblock Reaction-diffusion systems derived from kinetic theory for {M}ultiple
  {S}clerosis.
\newblock {\em Math. Models Methods Appl. Sci.}, 34(07):1279--1308, 2024.

\bibitem{othmer1988models}
H.~G. Othmer, S.~R. Dunbar, and W.~Alt.
\newblock Models of dispersal in biological systems.
\newblock {\em J. Math. Biology}, 26(3):263--298, 1988.

\bibitem{othmer2000diffusion}
H.~G. Othmer and T.~Hillen.
\newblock The diffusion limit of transport equations derived from velocity-jump
  processes.
\newblock {\em SIAM J. Appl. Math.}, 61(3):751--775, 2000.

\bibitem{othmer2002diffusion}
H.~G. Othmer and T.~Hillen.
\newblock The diffusion limit of transport equations {II}: {C}hemotaxis
  equations.
\newblock {\em SIAM J. Appl. Math.}, 62(4):1222--1250, 2002.

\bibitem{perez2012stochastic}
J.~P{\'e}rez-Vel{\'a}zquez, R.~Schlicht, G.~Dulla, B.~A. Hense, C.~Kuttler, and
  S.~E. Lindow.
\newblock Stochastic modeling of {P}seudomonas syringae growth in the
  phyllosphere.
\newblock {\em Math. Biosci.}, 239(1):106--116, 2012.

\bibitem{pusey2011antibiosis}
P.~Pusey, V.~Stockwell, C.~Reardon, and B.~Smits, T.H.M .and~Duffy.
\newblock Antibiosis activity of pantoea agglomerans biocontrol strain
  {E}{3}{2}{5} against {E}rwinia amylovora on apple flower stigmas.
\newblock {\em Phytopathology}, 101(10):1234--1241, 2011.

\bibitem{remus2012variation}
M.~N. Remus-Emsermann, R.~Tecon, G.~A. Kowalchuk, and J.~H. Leveau.
\newblock Variation in local carrying capacity and the individual fate of
  bacterial colonizers in the phyllosphere.
\newblock {\em ISME J.}, 6(4):756--765, 2012.

\bibitem{rionero2019hopf}
S.~Rionero.
\newblock Hopf bifurcations in dynamical systems.
\newblock {\em Ric. Mat.}, 68:811--840, 2019.

\bibitem{schlechter2019driving}
R.~O. Schlechter, M.~Miebach, and M.~N. Remus-Emsermann.
\newblock Driving factors of epiphytic bacterial communities: a review.
\newblock {\em Journal of Advanced Research}, 19:57--65, 2019.

\bibitem{schreiber2005plant}
L.~Schreiber, U.~Krimm, D.~Knoll, M.~Sayed, G.~Auling, and R.~M. Kroppenstedt.
\newblock Plant-microbe interactions: identification of epiphytic bacteria and
  their ability to alter leaf surface permeability.
\newblock {\em New Phytol.}, 166(2):589--594, 2005.

\bibitem{tecon2018cell}
R.~Tecon, A.~Ebrahimi, H.~Kleyer, S.~Erev~Levi, and D.~Or.
\newblock Cell-to-cell bacterial interactions promoted by drier conditions on
  soil surfaces.
\newblock {\em Proc. Nat. Acad. Sci.}, 115(39):9791--9796, 2018.

\bibitem{tulumello2014cross}
E.~Tulumello, M.~C. Lombardo, and M.~Sammartino.
\newblock Cross-diffusion driven instability in a predator-prey system with
  cross-diffusion.
\newblock {\em Acta Appl. Math.}, 132(1):621--633, 2014.

\bibitem{turing1990chemical}
A.~M. Turing.
\newblock The chemical basis of morphogenesis.
\newblock {\em Bull. Math. Biology}, 52:153--197, 1990.

\bibitem{van2013explaining}
A.~van~der Wal, R.~Tecon, J.~Kreft, W.~M. Mooij, and J.~H.~J. Leveau.
\newblock Explaining bacterial dispersion on leaf surfaces with an
  individual-based model ({PHYLLOSIM}).
\newblock {\em PLOS ONE}, 8(10):e75633, 2013.

\bibitem{wangersky1978lotka}
P.~J. Wangersky.
\newblock {L}otka-{V}olterra population models.
\newblock {\em Annu. Rev. Ecol. Evol. Syst.}, 9(1):189--218, 1978.

\bibitem{zagour2019modeling}
M.~Zagour.
\newblock {\em Modeling and mathematical analysis of complex systems: {K}inetic
  and macroscopic approaches and applications in biology and vehicular
  traffic}.
\newblock PhD thesis, Universit{\'e} Cadi Ayyad Marrakech (Maroc), 2019.

\bibitem{zhigun2022flux}
A.~Zhigun.
\newblock Flux limitation mechanisms arising in multiscale modelling of cancer
  invasion.
\newblock In {\em Mathematical Proceedings of the Royal Irish Academy}, volume
  122, pages 5--26. Royal Irish Academy, 2022.

\bibitem{zhigun2022novel}
A.~Zhigun and C.~Surulescu.
\newblock A novel derivation of rigorous macroscopic limits from a micro-meso
  description of signal-triggered cell migration in fibrous environments.
\newblock {\em SIAM J. Appl. Math.}, 82(1):142--167, 2022.

\end{thebibliography}
	

\end{document}